\documentclass[letterpaper, 11 pt, onecolumn]{ieeeconf}
\IEEEoverridecommandlockouts  
\usepackage{tikz}
\usepackage{nicefrac}
\usepackage{graphicx}
\usepackage{amsmath}
\usepackage{amsfonts}
\usepackage{amssymb}
\usepackage{mathrsfs}
\usepackage{fix2col}
\usepackage{latexsym}
\usepackage[dvips]{epsfig}
\usepackage[font=footnotesize]{subfig}
\usepackage{hyperref}
\usepackage{graphics}
\usepackage{epstopdf}
\usepackage{enumerate}
\usepackage{epstopdf}
\usepackage{arydshln}
\usepackage{mathdots}
\usepackage{wrapfig}
\usepackage{commands}
\usepackage{mathtools}
\DeclareMathAlphabet{\mathbbm}{U}{bbm}{m}{n}
\usepackage{algorithm,algcompatible}
\usepackage[normalem]{ulem} 

\newcommand{\vertiii}[1]{{\left\vert\kern-0.25ex\left\vert\kern-0.25ex\left\vert #1 
		\right\vert\kern-0.25ex\right\vert\kern-0.25ex\right\vert}}

\makeatletter
\newcommand{\opnorm}{\@ifstar\@opnorms\@opnorm}
\newcommand{\@opnorms}[1]{%
	\left|\mkern-1.5mu\left|\mkern-1.5mu\left|
	#1
	\right|\mkern-1.5mu\right|\mkern-1.5mu\right|
}
\newcommand{\@opnorm}[2][]{%
	\mathopen{#1|\mkern-1.5mu#1|\mkern-1.5mu#1|}
	#2
	\mathclose{#1|\mkern-1.5mu#1|\mkern-1.5mu#1|}
}
\makeatother

\begin{document}
\bstctlcite{IEEEexample:BSTcontrol}
\title{From self-tuning regulators to reinforcement learning and back again}
\author{Nikolai Matni,$^a$ Alexandre Proutiere,$^{b,\ast}$ Anders Rantzer,$^{c,\ast}$ Stephen Tu$^d$
	\thanks{$^a$ N. Matni is with the Dept. of Electrical and Systems Engineering, University of Pennsylvania, Philadelphia, PA. $^b$ A. Proutiere is with the Division of Decision and Control Systems, KTH Royal Institute of Technology, Stockholm, Sweden. $^c$ A. Rantzer is with the Department of Automatic Control, Lund University, Lund, Sweden. $^\ast$ A. Proutiere and A. Rantzer were partially supported by the Wallenberg AI, Autonomous Systems and Software Program (WASP) funded by the Knut and Alice
Wallenberg FoundationP.  $^d$ S. Tu is with Google Brain, New York City, NY.}}%
\maketitle
\begin{abstract}
Machine and reinforcement learning (RL) are increasingly being applied to plan and control the behavior of autonomous systems interacting with the physical world. Examples include self-driving vehicles, distributed sensor networks, and agile robots. However, when machine learning is to be applied in these new settings, the algorithms had better come with the same type of reliability, robustness, and safety bounds that are hallmarks of control theory, or failures could be catastrophic. Thus, as learning algorithms are increasingly and more aggressively deployed in safety critical settings, it is imperative that control theorists join the conversation. The goal of this tutorial paper is to provide a starting point for control theorists wishing to work on learning related problems, by covering recent advances bridging learning and control theory, and by placing these results within an appropriate historical context of system identification and adaptive control.
\end{abstract}

\section{Introduction \& Motivation}
\label{sec:intro}

With their recent successes in image classification, video game playing  \cite{atari-nature}, sophisticated robotic simulations \cite{duan2016benchmarking,Levine16}, and complex strategy games such as Go \cite{alphago,silver2018general}, machine and reinforcement learning (RL) are now being applied to plan and control the behavior of autonomous systems that interact with physical environments. Such systems, which include self-driving vehicles and agile robots, must interact with complex environments that are ever changing and difficult to model, strongly motivating the use of data-driven techniques. However, if machine learning is to be applied in these new settings, the resulting algorithms must come with the reliability, robustness, and safety guarantees that typically accompany results in the control theory literature, as failures could be catastrophic.  Thus, as RL algorithms are increasingly and more aggressively deployed in safety critical settings, control theorists must be part of the conversation.

To that end, it is important to recognize that while the applications areas and technical tools are new, the challenges faced -- uncertain and time varying systems and environments, unreliable sensing modalities, the need for robust stability and performance, etc. -- are not, and that many classical results from the system identification and adaptive control literature can be brought to bear on these problems.  In the case of discrete time linear systems, adaptive control algorithms further come with strong guarantees of asymptotic consistency, stability, and optimality, and similarly elucidate some of the fundamental challenges that are still being wrestled with today, such as rapidly identifying a system model (exploration) while robustly/optimally controlling it (exploitation).
 
Indeed, at a cursory glance, classical self-tuning regulators have the same objective as contemporary RL: an initial control policy and/or model is posited, data is collected, and a refined model/policy is produced, often in an online fashion. However, until recently, there has been relatively little contact between the two research communities. As a result, the tools and analysis objectives are different. One such feature, which will be the focus of this tutorial, is that 
RL and online learning algorithms are often analyzed in terms of \emph{finite-data} guarantees, and as such, are able to provide \emph{anytime guarantees} on the quality of the current behavior. Such finite-data guarantees are obtained by integrating tools from optimal control, stochastic optimization, and high-dimensional statistics -- whereas the first two tools are familiar to the controls community, the latter is less so.  A major theme will be that of uncertainty quantification: indeed the importance of relating uncertainty quantification to control objectives was emphasized already in the 1960-70s \cite{aastrom1971system}. Moreover, the fragility of certainty equivalent control was one of the motivating factors for the development of robust adaptive control methodologies. 

In this tutorial paper and our companion paper \cite{sysID}, we highlight recent advances that provide non-asymptotic analysis of adaptive algorithms.  Our aim is for these papers is for them to serve as a jumping off point for control theorists wanting to work in RL problems.  In \cite{sysID}, we present an overview of tools and results on finite-data guarantees for system identification.  This paper focuses on finite-data guarantees for self-tuning and adaptive control strategies, and is structured as follows:
\begin{itemize}
\item \textbf{Section \ref{sec:lit-review}}: provides an extensive literature review of work spanning classical and recent results in system identification, adaptive control, and RL.
\item \textbf{Section \ref{sec:fundamentals}}: introduces the fundamental problem and performance metrics considered in RL, and relates them to examples familiar to the controls community.
\item \textbf{Section \ref{sec:mdps}}: provides a survey of contemporary results for problems with finite state and action spaces.
\item \textbf{Section \ref{sec:model_based_lqr}}: shows how system estimates and error bounds can be incorporated into model-based self-tuning regulators with finite-time performance guarantees.
\item \textbf{Section \ref{sec:model_free}}: presents guarantees for model-free methods, and shows that a complexity gap exists between model-based and model-free methods.
\end{itemize}

\section{Literature Review}
\label{sec:lit-review}

The results we present in this paper draw heavily from three broad areas of control and learning theory: system identification, adaptive control, and approximate dynamic programming (ADP) or, as it has come to be known, reinforcement learning. Each of these areas has a long and rich history and a general literature review is outside the scope of this tutorial. Below we will instead emphasize pointers to good textbooks and survey papers, before giving a more careful account of recent work.

\subsubsection{System Identification}
The estimation of system behavior from input/output experiments has a well-developed theory dating back to the 1960s, particularly in the case of linear-time-invariant systems. Standard reference texts on the topic include \cite{aastrom1971system,ljung1999system,chen2012identification,goodwin1977dynamic}. The success of discrete time series analysis by Box and Jenkins \cite{box2015time} provided an early impetus for the extension of these methods to the controlled system setting. 
Important connections to information theory were established by Akaike \cite{akaike1973information}. 
 The rise of robust control in the 1980s further inspired system identification procedures, wherein model errors were optimized
 under the assumption of adversarial noise processes \cite{makila1995worst}. 
Another important step was the development of subspace methods \cite{van2012subspace}, which became a powerful tool for identification of multi-input multi-output systems. 

\subsubsection{Adaptive Control}
System identification on real systems can be tedious, time consuming, and require skilled personnel. Adaptive control could then offer a simpler path forward. Moreover, adaptation offers an effective way to compensate for time-variations in the system dynamics. An early driving application was aircraft autopilot development in the 1950s. Aerospace applications needed control strategies that automatically compensated for changes in dynamics due to altitude, speed, and flight configuration  \cite{stein1980adaptive}. Another important early application was ship steering, where adaptation is used to compensate for wave effects \cite{aastrom1980use}. Standard textbooks include \cite{aastrom2013adaptive,goodwin2014adaptive,narendra2012stable,sastry2011adaptive,astolfi2007nonlinear}. 

A direct approach to adaptive control expounded by Bellman \cite{bellman2015adaptive} was to tackle the problem using dynamic programing by augmenting the state to contain the conditional distribution of the unknown parameters.  From this it was seen that in such problems, control served the dual purpose of exciting the system to aid in its identification -- hence the term dual control.  This approach suffered from the curse of dimensionality, which lead to the development of approximation techniques \cite{bertsekas1996neuro,powell2007approximate} that ultimately evolved into modern day RL.  

A more practical approach was self-tuning adaptive control, pioneered by \cite{kalman1958design,aastrom1973self} and followed by a long sequence of contributions to adaptive control theory, deriving conditions for convergence, stability, robustness and performance under various assumptions. For example, \cite{ljung77ana} analysed adaptive algorithms using averaging, and \cite{goodwin1981discrete} derived an algorithm that gives mean square stability with probability one. On the other hand, conditions that may cause instability were studied in \cite{ega79book}, \cite{ioannou1984instability} and \cite{rohrs1985robustness}. Finally, \cite{guo1995convergence} gave conditions for optimal asymptotic rates of convergence.  More recent adaptive approaches include the L1 adaptive controller \cite{hovakimyan2010L1}, and model free adaptive control \cite{hou2011data,hou2013model,fliess2013model}.

\subsubsection{Automatic Tuning and Repeated Experiments}
There is also an extensive literature on automatic procedures for initialization (tuning) of controllers, without further adaptation after the tuning phase. A successful example is 
auto-tuning of PID controllers \cite{aastrom1984automatic}, where a relay provides non-linear feedback during the tuning phase. Another important tuning approach, well established from an engineering perspective \cite{hjalmarsson1998iterative}, is based on repeated experiments with linear time-invariant controllers. Theoretical bounds on a related approach were obtained by Lai and Robbins \cite{lai1982iterated}. Specifically, they showed that a pseudo-regret of the state variance is lower bounded by $\Omega(\log(T))$.   Subsequent work by Lai \cite{lai1986asymptotic} showed that this bound was tight.  Recently, Raginsky \cite{raginsky2010divergence} revisited this problem formulation, and showed that for any persistently exciting controller, the time taken to achieve state variance less than $\epsilon$ is at least $\Omega(\tfrac{n^2}{\epsilon}\log(1/\epsilon))$ for a system of dimension~$n$.



\subsubsection{Dynamic Programming and Reinforcement Learning}
A major part of the literature on dynamic programming is devoted to ``tabular MDPs,'' i.e. 
systems for which the state and action spaces are discrete and small enough to be stored in memory. The classic texts \cite{sutton2018reinforcement,bertsekas1996neuro,powell2007approximate} highlight computationally efficient approximation techniques for solving these problems.  They include  Monte Carlo methods, temporal-difference (TD) learning \cite{tesauro1995temporal,sutton1988learning} (which encompass SARSA \cite{rummery1994line} and Q-learning \cite{watkins1989learning,watkins1992q,barnard1993temporal}), value and Q-function approximation via Neural Networks, kernel methods, least-squares TD (LSTD) \cite{barto1995learning,barnard1993temporal,nedic2003least}, and policy gradient methods such as REINFORCE \cite{ watkins1992q,phansalkar1995local}  and Actor-Critic Methods \cite{barto1983neuronlike,degris2012off}.  

Recent advances in both algorithms and computational power have allowed RL methods to solve incredibly complex tasks in very large discrete spaces that far exceed the tabular setting, including video games \cite{atari-nature}, Go \cite{alphago}, chess, and shogi \cite{silver2018general}.  This success has renewed an interest in applying traditional model-free RL methods, such as Q-learning \cite{lillicrap2015continuous} and policy optimization \cite{schulman2015trust}, to continuous problems in robotics \cite{duan2016benchmarking,Levine16}.  Thus far, however, the deployment of systems trained in this way has been limited to simulation environments \cite{todorov2012mujoco} or highly controlled laboratory settings, as the training process for these systems is both data hungry and highly variable \cite{mania2018simple}.
%

\subsubsection{System Identification Revisited}

With few exceptions (e.g., \cite{dahleh1993sample}), prior to the 2000s, the literature on system identification and adaptive control focused on asymptotic error characterization and consistency guarantees. In contrast, contemporary results in statistical learning seek to characterize \emph{finite time and finite data} rates, leaning heavily on tools from stochastic optimization and concentration of measure.  Such finite-time guarantees provide estimates of both system parameters and their uncertainty, allowing for a natural bridge to robust/optimal control. Early such results, characterizing rates for parameter identification \cite{campi2002finite,vidyasagar2006learning}, featured conservative bounds which are exponential in the system degree and other relevant quantities. More recent results, focused on state-space parameter identification for LTI systems, have significantly improved upon these bounds.  In \cite{hardt2018gradient}, the first polynomial time guarantees for identifying a stable linear system were provided -- however, these guarantees are in terms of predictive output performance of the model, and require rather stringent assumptions on the true system. In \cite{learning-lqr}, it was shown, assuming that the state is directly measurable and the system is driven by white in time Gaussian noise, that solving a least-squares problem using independent data points taken from different trials achieves order optimal rates that are linear in the system dimension. This result was generalized to the single trajectory setting for (i) marginally stable systems in \cite{simchowitz2018learning}, (ii) unstable systems  in \cite{sarkar2018fast},
and (iii) partially observed stable systems in \cite{oymak2018non,sarkar2019finite,tsiamis2019finite}.  We note that analogous results also exist in the fully observed setting for the identification of sparse state-space parameters \cite{fattahi2018data,fattahi2019learning}, where rates are shown to be logarithmic in the ambient dimension, and polynomial in the number of nonzero elements to be estimated.

\subsubsection{Automatic Tuning Revisited}

There has been renewed interest, motivated in part by the expansion of reinforcement learning to continuous control problems, in the study of automatic tuning as applied to the Linear Quadratic Regulator.  More closely akin to iterative learning control \cite{bristow2006survey}, Fietcher \cite{fiechter1997pac} showed that the discounted LQR problem is Probably Approximately Correct (PAC) learnable in an episodic setting.  In \cite{learning-lqr}, Dean et al. dramatically improved the generality and sharpness of this result by extending it to the traditional infinite horizon setting, and leveraging contemporary tools from concentration of measure and robust control.  

\subsubsection{Adaptive Control Revisited}
Contemporary results tend to draw on ideas from the bandits literature.  A non-asymptotic study of the adaptive LQR problem was initiated by Abbasi-Yadkori and Szepesvari~\cite{abbasi2011regret}.  They use an Optimism in the Face of Uncertainty (OFU) based approach,
where they maintain confidence ellipsoids of system parameters and select those parameters that lead to the best closed loop performance. 
{While the OFU method achieves the optimal $O(T^{1/2})$ regret, solving the OFU sub-problem is computationally challenging.}  To address this issue, other exploration methods were studied.  Thompson sampling \cite{russo17} is used to achieve $O(T^{1/2})$ regret for scalar systems \cite{abeille18}, and \cite{ouyang17} studies a Bayesian setting with a particular Gaussian prior.  
%
%
Both \cite{abbasi2019model} and \cite{dean2018regret} give tractable algorithms
which achieve sub-linear frequentist regret of $O(T^{2/3})$ without the Bayesian setting of \cite{ouyang17}.
Follow up work \cite{cohen2019learning} showed that this rate could be improved to $O(T^{1/2})$ by leveraging a novel semi-definite relaxation.  More recently, \cite{mania2019certainty} show that as long as the initial system parameter estimates are sufficiently accurate, certainty equivalent (CE) control achieves $O(T^{1/2})$ regret with high probability. Finally, Rantzer \cite{rantzer2018concentration} shows that for a scalar system with a minimum variance cost criterion, a simple self-tuning regulator scheme achieves $O(\log(T))$ expected regret after an initial burn in period, thus matching the lower bound established by \cite{lai1982iterated}.  

Much of the recent work addressing the sample complexity of the LQR problem was motivated by the desire to understand RL algorithms on a simple baseline \cite{recht2018tour}.  In addition to the model-based approaches described above, model-free methods have also been studied.  Model-free methods for the LQR problem were put on solid theoretical footing in \cite{bradtke1994adaptive}, where it was shown that controllability and persistence of excitation were sufficient to guarantee convergence to an optimal policy. In \cite{tu2018least}, the first finite time analysis for LSTD as applied to the LQR problem is given, in which they show that $O(n^3/\epsilon^2)$ samples are sufficient to estimate the value function up to $\epsilon$-accuracy, for $n$ the state dimension. 
%
%
Subsequently, Fazel et al.~\cite{fazel18} also showed that randomized search algorithms similar to policy gradient
can learn the optimal controller with a polynomial number of samples in the noiseless case; however
an explicit characterization of the dependence of the sample complexity on the parameters
of the true system was not given, and the algorithm is dependent on the knowledge of an initially stabilizing controller.
Similarly, Malik et al.~\cite{malik2019dfo} study the behavior of random finite differencing for LQR.
Finally, \cite{tu2018gap} shows that there exists a family of systems for which there is a sample complexity gap of at least a factor of state dimension between LSTD/policy gradient methods and simple CE model-based approaches.


\section{Fundamentals}  
\label{sec:fundamentals}
We study the behavior of \emph{Markov Decision Processes (MDP)}.  In the finite horizon setting of length $T$, we consider 
\begin{equation}
\begin{array}{rl}
\min_{\tf \pi} & \E\left[\sum_{t=0}^{T-1} c_t(x_t,u_t) + c_T(x_T)\right] \\
\mathrm{s.t.} & x_{t+1} = f_t(x_t, u_t, w_t) \\
& u_t = \pi_t(x_{0:t},u_{0:t-1}), 
\end{array}
\label{eq:MDP-T}
\end{equation}
for $x_t \in \mathcal{X}$ the system state, $u_t \in \mathcal{U}$ the control input, $w_t \in \mathcal{W}$ the state transition randomness, and $\tf \pi = \{\pi_0, \pi_1, \dots, \pi_{T-1}\}$ the control policy, with $\pi_t:\mathcal{X}^t \times \mathcal{U}^{t-1} \to \mathcal{U}$ a possibly random mapping.  With slight abuse of notation, we will use $n_x$, $n_u$, and $n_w$ to denote (i) the dimension of $\mathcal{X}$, $\mathcal{U}$, and $\mathcal{W}$, respectively, when considering continuous state and action spaces, and (ii) the cardinality of  $\mathcal{X}$, $\mathcal{U}$, and $\mathcal{W}$, respectively, when considering discrete state and action spaces.

We consider settings where both the cost functions $\{c_t\}_{t=0}^T$ and the dynamics functions $\{f_t\}_{t=0}^{T}$ may not be known.  Finally, we assume that the primitive random variables $(x_0, w_{0:T})$ are defined over a common probability space with known and independent distributions -- the expectation in the cost is taken with respect to these and the policy $\tf \pi$.

\subsection{Dynamic Programming Solutions}
When the transition functions $\{f_t\}$ and costs $\{c_t\}$ are known, problem \eqref{eq:MDP-T} can be solved using dynamic programming.  As the dynamics are Markovian, we restrict our search to policies of the form $u_t = \pi_t(x_t)$ without loss of optimality.  Define the value function of problem \eqref{eq:MDP-T} at time $t$ to be
\begin{equation}
\begin{array}{rcl}
V_T(x_T) &=& \E[c_T(x_T)]\\
V_{t}(x_t) &=& \min_{u_t}\E\left[c_t(x_t,u_t) + V_{t+1}(f_t(x_t,u_t,w_t))\right].
\end{array}
\label{eq:finite-value-fcn}
\end{equation}
Iterating through this process yields both an optimal policty $\tf \pi^\star$, and the optimal cost-to-go $V_0(x_0)$ that it achieves.

Moving to the infinite horizon setting, we assume that the cost function and dynamics are static, i.e., that $c_t(x_t,u_t) \equiv c(x_t,u_t)$ and $f_t\equiv f(x_t,u_t,w_t)$ for all $t \geq 0$.

We begin by introducing the discounted cost setting, wherein the cost-functional in optimization \eqref{eq:MDP-T} is replaced with
\begin{equation}
\E\left[\sum_{t=0}^{\infty} \gamma^t c(x_t,u_t) \right],
\label{eq:MDP-gamma}
\end{equation} 
for some $\gamma \in (0,1]$.  Note that if $c(x_t,u_t)$ is bounded almost surely and $\gamma < 1$, then the infinite sum \eqref{eq:MDP-gamma} is guaranteed to remain bounded, greatly simplifying analysis.  In this setting, evaluating the performance achieved by a fixed policy $\tf \pi$ consists of finding a solution to the following equation:
\begin{equation}
V_{\pi}(x) = \E\left[c(x,\pi(x)) + \gamma V_\pi(f(x,\pi(x),w))\right],
\label{eq:discounted-policy-evalutation}
\end{equation}
from which it follows immediately that the optimal value function $V_\star$ will satisfy
\begin{equation}
V_{\star}(x) = \min_u \E\left[c(x,u) + \gamma V_{\star}(f(x,u,w))\right],
\end{equation}
and $u_\star = \pi_\star(x)$.  One can show that under mild technical assumptions, iterative procedures such as policy iteration and value iteration {will converge to the optimal policy}.  

Next, we consider the asymptotic average cost setting, in which case the cost-functional in problem \eqref{eq:MDP-T} is set to 
\begin{equation}
\E\left[\lim_{T\to \infty} \frac{1}{T}\sum_{t=0}^{T-1} c(x_t,u_t) + c_T(x_T)\right].
\label{eq:MDP-inf}
\end{equation} 
Care must be taken to ensure that the limit converges, thus somewhat complicating the analysis -- however, this cost functional is often most appropriate for guaranteeing the stability for stochastic optimal control problems.

%

\begin{example}[Linear Quadratic Regulator]
Consider the Linear Quadratic Regulator (LQR), a classical instantiation of MDP \eqref{eq:MDP-T} from the optimal control literature, for $Q \succeq 0$ and $R \succ 0$:
\begin{equation}
\begin{array}{rl}
\min_{\tf \pi} & \frac{1}{T}\E\left[\sum_{t=0}^{T-1} x_t^\top Q x_t + u_t^\top R u_t + x_T^\top Q_T x_T \right] \\
\mathrm{s.t.} & x_{t+1} = Ax_t + Bu_t + w_t \\
& u_t = \pi_t(x_{0:t},u_{0:t-1}),
\end{array}
\label{eq:LQR-T}
\end{equation}
where the system state $x_t \in \R^{n_x}$, the control input $u_t \in \R^{n_u}$, and the disturbance process $w_t \in \R^{n_x}$ are independently and identically distributed as zero mean Gaussian random variables with a known covariance matrix $\Sigma_w$.

For a finite horizon $T$ and known matrices $(A,B)$, this problem can be solved directly via dynamic programming, leading to the optimal control policy
\begin{equation}
u_t^\star = -(B^\top P_{t+1} B + R)^{-1} B^\top P_{t+1} A x_t,
\label{eq:lqr-opt}
\end{equation}
where $P_t \succeq 0$ satisfies the Discrete Algebraic Riccati (DAR) Recursion initialized at $P_T = Q_T.$  Further, when the triple $(A,B,Q^{1/2})$ is stabilizable and detectable, the closed loop system is stable and hence converges to a stationary distribution, allowing us to consider the asymptotic average cost setting \eqref{eq:MDP-inf}, at which point the optimal control action is a static policy, defined as in \eqref{eq:lqr-opt}, but with $P_{t} \to P$, for $P\succeq0$ a solution of the corresponding DAR Equation.
\end{example}

\begin{example}[Tabular MDP]
Consider the setting where the state-space $\mathcal{X}$, the control space $\mathcal{U}$, and the disturbance process space $\mathcal{W}$ have finite cardinalities of $n_x$, $n_u$, and $n_w$, respectively, and further suppose that the underlying dynamics are governed by transition probabilities $P(x_{t+1} = x' | x_t, u_t)$.  We assume that the cardinalities $n_x$, $n_u$, and $n_w$ are such that $\mathcal{X}$, $\mathcal{U}$, and $\mathcal{W}$ can be stored in tabular form in memory and worked with directly.  These then induce the dynamics functions:
\begin{equation}
x_{t+1} = w_t(x_t,u_t),
\end{equation}
where $w_t(x_t,u_t) = x'$ with probability $P(x_{t+1} = x' | x_t, u_t)$.

In the case of average cost, for simplicity, we restrict our attention to {\it communicating} and {\it ergodic} MDPs. The former correspond to scenario where for any two states, there exists a stationary policy leading from one to the other with positive probability. For the latter, any stationary policy induces an ergodic Markov chain. In the average cost setting, one wishes to minimize $\lim_{T\to\infty} {1\over T}V_T(x)$. More precisely, the objective is to identify as fast as possible a stationary policy $\pi$ with maximal {\it gain} function $g^\pi$: for any $x\in {\cal X}$, $g^\pi(x):=\lim_{T\to\infty} {1\over T}V_T^\pi(x)$ where $V_T^\pi(x)$ denotes the average cost under $\pi$ starting in state $x$ over a time horizon $T$. When the MDP is ergodic, this gain does not depend on the initial state. To compute the gain of a policy, we need to introduce the {\it bias} function $h^\pi(x):=\textnormal{{\em C}-}\lim_{T\to\infty} \mathbb{E}^\pi[\sum_{t = 1}^\infty (c(x_t,u_t)  - g^\pi(x_t))|x_0=x]$ (where $\textnormal{{\em C}-}\lim_{T\to\infty}$ is the Cesaro limit) that quantifies the advantage of starting in state~$x$. $g^\pi$ and $h^\pi$ satisfy for any $x$:
$$
g^\pi(x)+h^\pi(x) = c(x,\pi(x)) +\sum_y p(y|x,\pi(x)) h^\pi(y).
$$    
The gain and bias functions $g^\star(x)$ and $h^\star(x)$ of an optimal policy verify Bellman's equation: for all $x$,
$$
g^\star(x)+h^\star(x) = \min_{u\in {\cal U}}\left( c(x,u) +\sum_y p(y|x,u) h^\star(y)\right).
$$
$h^\star$ is defined up to an additive constant.

When the transition probabilities and cost function are known, this problem can then be solved via value-iteration, policy-iteration, and linear programming.
\end{example}

\subsection{Learning to Control MDPs with Unknown Dynamics}
Thus far we have considered settings where the dynamics $\{f_t\}$ and costs $\{c_t\}$ are known.  Our main interest is understanding what should be done when these models are not known.  Our study will focus on the previous two examples, namely LQR and the tabular MDP setting.  While much of current work in reinforcement learning focusses on \emph{model-free} methods, we adopt a more control theoretic perspective on the problem and study model-based methods wherein we attempt to approximately learn the system model $\{f_t\}_{t=0}^T$, and then subsequently use this approximate  model for control design.  

Before continuing, we distinguish between episodic and single-trajectory settings.  An \emph{episodic task} is akin to traditional \emph{iterative learning control}, wherein a  task is repeated over a finite horizon, after which point the episode ends, and the system is reset to begin the next episode.  In contrast, a \emph{single-trajectory task} is akin to traditional \emph{adaptive control}, in that no such resets are allowed, and a single evolution of the system under an adaptive policy is studied.  

An underlying tension exists between identifying an unknown system and controlling it.  Indeed it is well known that without sufficient \emph{exploration} or \emph{excitation}, an incorrect model will be learned, possibly leading to suboptimal and even unstable system behavior; however, this exploration inevitably degrades system performance.  Informally, this tension leads to a fundamental tradeoff between how quickly a model can be learned, and how well it can be controlled during this process.  Current efforts seek to explicitly address and quantify these tradeoffs through the use of performance metrics such as the Probably Approximately Correct (PAC) and Regret frameworks, which we define next.  For episodic tasks, we assume the horizon of each episode to be of length $H$, and consider guarantees on performance as a function of the number of episodes $T$ that have been evaluated.  For single trajectory tasks, we consider infinite horizon problems, and the definitions provided are equally applicable to the discounted and asymptotic average cost settings.  The definitions that follow are adapted from \cite{fiechter1997pac, dann2017unifying, strehl2006pac}, among others.




\subsection{PAC-Bounds}
\paragraph{Episodic PAC-Bounds}
We consider episodic tasks over a horizon $H$, where $H$ may be infinite but the user is allowed to reset the system at a prescribed time $H_r$.  Let $V_{\star}$ be the optimal cost achievable, and $N_\epsilon$ be the number of episodes for which $\tf \pi$ is not $\epsilon$-optimal, i.e., the number of episodes for which $V_{\tf\pi}  >  V_\star + \epsilon$.  Then, a policy $\tf \pi$ is said to be episodic-$(\epsilon,\delta)$-PAC if, after $T$ episodes, it satisfies\footnote{We note that in the discounted setting, we also ask that the bound on $N_\epsilon$ depend polynomially on $1/(1-\gamma)$.  We also note that modern definitions of PAC-learning \cite{dann2017unifying} require that the bound on $N_\epsilon$ depend polynomially on $\log(1/\delta)$ -- this is a reflection of results from contemporary high-dimensional statistics that allow for more refined concentration of measure guarantees.  Finally, the polynomial dependence on the horizon $H$ is only enforced for finite horizons $H$ -- in the case of an infinite horizon task, $N_\epsilon$ must not depend on the $H$.}
\begin{equation}
\Prob{N_\epsilon > \mathrm{poly}(n_x, n_u, H, 1/{\epsilon},1/\delta)}\leq \delta.
\end{equation}

These guarantees state that the chosen policy is $\epsilon$-optimal on all but a  number of episodes polynomial in the problem parameters, with probability at least $1-\delta$.  Many $(\epsilon,\delta)$ PAC algorithms operate in two phases: the first is solely one of exploration so as to identify an approximate system model, and the second is solely one of exploitation, wherein the approximate system model is used to synthesize a control policy.  Therefore, informally one can view PAC guarantees as characterizing the number of episodes needed to identify a model that can be used to synthesize an $\epsilon$-optimal policy.


\begin{example}[LQR is episodic PAC-Learnable]
The results in \cite{learning-lqr} imply that the LQR problem with an asymptotic  average cost is episodic PAC-learnable.  In particular, it was shown that a simple open-loop exploration process of injecting white in time Gaussian noise over at most $\mathrm{poly}(n_x,n_u,H_r,1/\epsilon,\log(1/\delta))$ episodes, followed by a least-squares system identification and uncertainty quantification step, can be used with a robust synthesis method to generate a policy $\tf \pi$  which guarantees that
\begin{equation}
\Prob{ V_{\tf \pi} - V_\star \geq \epsilon} \leq \delta,
\end{equation}
when the LQR problem is initialized at $x_0 = 0$.  Hence the resulting algorithm meets the modern definition of being $(\epsilon,\delta)$-PAC-learnable.  We revisit this example in Section \ref{sec:model_based_lqr}.
\end{example}

\paragraph{Single-trajectory PAC-Bounds}
We consider single-trajectory tasks over an infinite horizon, and let $V_{\tf \pi}(x_t)$ denote the cost-to-go from state $x_t$ achieved by a policy $\tf \pi$, and $V_{\star}(x_t)$ be the optimal cost-to-go achievable. We further let $N_\epsilon$ be the number of time-steps for which $\tf \pi$ is not $\epsilon$-optimal in either an absolute or relative sense, i.e., the number of time-steps for which $V_{\tf\pi}(x_t)  >  V_\star(x_t) + \epsilon$ or  $V_{\tf\pi}(x_t)  >  V_\star(x_t)(1 + \epsilon)$, respectively.  Then, a policy $\tf \pi$ is said to be $(\epsilon,\delta)$-PAC if it satisfies\footnote{We make the same modifications to this definition for the discounted case and the dependence on $1/\delta$ as in the episodic setting.}
\begin{equation}
\Prob{N_\epsilon > \mathrm{poly}(n_x, n_u, 1/{\epsilon},1/\delta)}\leq \delta.
\end{equation}
These guarantees should be interpreted as saying that the chosen policy is at worst $\epsilon$-suboptimal on all but $\mathrm{poly}(\tfrac{1}{\epsilon},\log(1/\delta))$ time-steps, with probability at least $1-\delta$.  As in the episodic setting, one can view these PAC guarantees as characterizing the number of time-steps needed to identify a model that can be used to synthesize an $\epsilon$-optimal policy.



\paragraph{Limitations of PAC-Bounds} As an algorithm that is $(\epsilon,\delta)$-PAC is only penalized for suboptimal behavior exceeding the $\epsilon$ threshold, there is no guarantee of convergence to an optimal policy.  In fact, as pointed out in \cite{dann2017unifying} and illustrated in the LQR example above, many PAC algorithms cease learning once they are able to produce an $\epsilon$-suboptimal strategy.

\subsection{Regret Bounds}
We focus on regret bounds for the single-trajectory setting, as this is the most common type of guarantee found in the literature, but note that analogous episodic definitions exist (cf., \cite{dann2017unifying}).  The regret framework evaluates the quality of an adaptive policy by comparing its running cost to a suitable baseline.  Let $b_T$ represent the baseline cost at time $T$, and define the regret incurred by a policy $\tf\pi = \{\pi_0, \pi_1, \dots \}$ to be
\begin{equation}
R^\pi(T):= \sum_{t=0}^{T} c_t(x_t,\pi_t(x_{0:t},u_{0:t-1})) - b_T.
\label{eq:regret}
\end{equation}
Note that $b_T$ is user specified, and is often chosen to be the expected optimal cost achievable by a policy with full knowledge of the system dynamics. The two most common regret guarantees found in the literature are expected regret bounds, and high probability regret bounds.  In the expected regret setting, the goal is to show that 
\begin{equation}
\mathbb E R^\pi(T) \leq \mathrm{poly}(n_x,n_u,T),
\end{equation}
whereas in the high-probability regret setting, the goal is to show that\footnote{As in the PAC setting, modern definitions often require the dependence to be polynomial in $\log(1/\delta)$.}
\begin{equation}
\Prob{R^\pi(T) \geq \mathrm{poly}(n_x,n_u,T,1/\delta)} \leq \delta.
\end{equation}
These bounds therefore quantify the rate of convergence of the cost achieved by the adaptive policy to the baseline cost, providing \emph{any time guarantees} on performance relative to a desirable baseline.  From the definition of $R^\pi(T)$, it is clear that one should strive for an $o(T)$ dependence, as this implies that the cost achieved by the adaptive policy converges with at least sub-linear rate to the base cost $b_T$.  Further,  in contrast to the PAC framework, \emph{all sub-optimal behavior} is tallied by the running regret sum, and hence exploration and exploitation must be suitably balanced to achieve favorable bounds.



\begin{example}[Regret bounds for LQR]
The study of regret bounds for LQR was initiated in \cite{abbasi2011regret}.  Here we summarize a recent treatment of the problem, as provided in \cite{mania2019certainty}.  There, the authors study the performance of CE control for LQR, and study a regret measure of the form
\begin{equation}
R^\pi(T) := \sum_{t=0}^T x_t^\top Q x_t + u_t^\top R u_t - TV_\star,
\end{equation}
with $V_\star := \min_u \E\left[\lim_{T\to \infty} T^{-1}\sum_{t=0}^T x_t^\top Q x_t + u_t^\top R u_t\right]$ the optimal asymptotic average cost achieved by the true optimal LQR controller.  They show that  the control policy $u_t = \hat{K}x_t + \eta_t,$
which has an exploration term $\eta_t \sim \Normal (0,\sigma_{\eta,t}^2 I)$ added to the CE controller, achieves  the regret bound
\begin{equation}
R^\pi(T) \leq \mathrm{poly}(n_x,n_u,\log(1/\delta))O(T^{1/2}),
\end{equation}
with probability at least $1 - \delta$ so long as $\sigma_{\eta,t}^2 \sim t^{-1/2}$ and the initial estimates of the system dynamics $(\Ahat,\Bhat)$ are sufficiently accurate.  We revisit this example in Section \ref{sec:model_based_lqr}.
\end{example}

\paragraph{Limitations of Regret Bounds} As regret only tracks the integral of suboptimal behavior, it does not distinguish between a few severe mistakes and many small ones. In fact, \cite{dann2017unifying} shows that for Tabular MDP problems, an algorithm achieving optimal regret may still make infinitely many mistakes that are maximally suboptimal.  Thus regret bounds cannot provide guarantees about transient worst-case deviations from the baseline cost $b_T$, which may have implications on guaranteeing the robustness or safety of an algorithm.  We comment further on regret for discrete MDPs in the next section.


\section{Optimal control of unknown discrete systems}\label{sec:mdps}
\label{sec:mdps}

This section addresses reinforcement learning for stationary MDPs with finite state and control spaces of respective cardinalities $n_x$ and $n_u$. When the system is in state $x$ and the control input is $u$, the system evolves to state $x'$ with probability $p(x'|x,u)$, and the cost $c_t(x,u)$ induced is independently drawn from a distribution $q(\cdot |x,u)$ with expectation $c(x,u)$. Costs are bounded, and for any $(x,u)$, the distribution $q(\cdot|x,u)$ is absolutely continuous w.r.t. to a measure $\lambda$ (for example, Lebesgue measure). 

We consider both average-cost and discounted MDPs (refer to Example 2), and describe methods (i) to derive fundamental performance limits (in terms of regret and sample complexity), and (ii) to devise efficient learning algorithms. The way the learner samples the MDP may significantly differs in the literature, depending on the objective (average or discounted cost), and on whether one wishes to derive fundamental performance limits or performance guarantees of a given algorithm. For example, in the case of average cost, typically, the learner gathers information about the system in an online manner following the system trajectory, a sampling model referred previously to as the single trajectory model. Most sample complexity ananlyses are on the contrary derived under the so-called {\it generative} model, where any state-control pair can be sampled in $O(1)$ time. Generative models are easier to analyze but hide the difficult issue of navigating the state space to explore various state-control pairs.   

\subsection{Average-cost MDPs}

For average-cost MDPs, we are primarily interested in devising algorithms with minimum regret, defined for a given learning algorithm $\pi$ as
$$
R^\pi(T) :=  \sum_{t=1}^T c_t(x_t^\pi,u_t^\pi) -\sum_{t=1}^T c_t(x_t^\star,u_t^\star),
$$
where $x_t^\pi$ and $u_t^\pi$ are the state and control input under the policy $\pi$ at time-step $t$, and similarly the superscript $\star$ corresponds to an optimal stationary control policy. Next, we discuss the type of regret guarantees we could aim at.

\medskip
\noindent
{\it Expected vs. high-probability regret.} We would ideally wish to characterize the complete regret distribution. This is however hard, and most existing results provide guarantees either in expectation or with high probability. In the case of finite state and control spaces, guarantees in high-probability can be easy to derive and not very insightful. Consider for example, a stochastic bandit problem (a stateless MDP) where the goal is to identify the control input with the lowest expected cost. An algorithm exploring each control input at least $\log(1/\delta)$ times would yield a regret in ${\cal O}(\log(1/\delta))$ with probability greater than $1-\delta$ (this is a direct application of Hoeffding's inequality). This simple observation only holds for a fixed MDP (the gap between the costs of the various inputs cannot depend on $\delta$ nor on the time at which regret is evaluated), and relies on the assumption of bounded costs. Even in the simplistic stochastic bandit problem, analyzing the distribution of the regret remains an open and important challenge, refer to \cite{Audibert2009, Salomon2011} for initial discussions and results.

\medskip
\noindent
{\it Problem-specific vs. minimax regret guarantees.} A regret upper bound is problem-specific if it explicitly depends on the parameters defining the MDP. Such performance guarantees capture and quantify the hardness of learning to control the system. Minimax regret guarantees are far less precise and informative, since they concern the {\it worst} system among possibly all systems. An algorithm with good minimax regret upper bound behaves well in the worst case, but does not necessarily learn and adapt to the system it aims at controlling. 

\medskip
Guided by the above observations, we focus on the expected regret, and always aim, when this is possible, at deriving problem-specific performance guarantees. 

\medskip

\subsubsection{Regret lower bounds} 

We present here a unified simple method to derive both problem-specific and minimax regret lower bounds. This method has been developed mainly in the bandit optimization literature \cite{kaufmann2016complexity,CombesP14b,Combes2014,garivier193} as a simplified alternative to Lai and Robbins techniques \cite{lai1985}.    

Let $\phi=(p,q)$ denote the true MDP. Consider a second MDP $\psi=(p',q')$. For a given learning algorithm $\pi$, define by ${\cal L}^\pi(T)$ the log-likelihood ratio of the corresponding observations under $\phi$ and $\psi$. By a simple extension of the Wald's lemma, we get:
\begin{equation}\label{eq:low1}
\mathbb{E}_\phi^\pi[{\cal L}^\pi(T)] = \sum_{x,u} \mathbb{E}_{\phi}^\pi[N_{xu}(T)]KL_{\phi|\psi}(x,u),
\end{equation}
where $N_{xu}^\pi(T)$ is the number of times the state-control pair $(x,u)$ is observed, where $\mathbb{E}_\phi^\pi$ is the expectation taken w.r.t. to the distrbution of observations made under $\pi$ for the MDP $\phi$, and where $KL_{\phi |\psi}(x,u)$ is the KL divergence between the distributions of the observations made in $(x,u)$ under $\phi$ and $\psi$. These observations concern the next state, and the realized reward, and hence:
\begin{align*}
KL_{\phi|\psi}(x,u) = \sum_{y}& p(y | x,u)  \log \frac{p(y | x,u) }{p'(y | x,u) } \\&\quad + \int q(r | x,u) \log \frac{q(r | x,u)}{q'(r | x,u) }\lambda(dr).
\end{align*}
Now the data processing inequality states that for any event $E$ or any $[0,1]$-valued random variable depending on all observations up to time $T$, i.e., ${\cal F}_T^\pi$-measurable (${\cal F}_T^\pi$ is the $\sigma$-algebra generated by the observations under $\pi$ up to time $T$):
$$
\mathbb{E}_\phi^\pi[{\cal L}^\pi(T)]  \ge \max\{ kl(\mathbb{P}_\phi^\pi[E],  \mathbb{P}_\psi^\pi[E]), kl(\mathbb{E}_\phi^\pi[Z],  \mathbb{E}_\psi^\pi[Z])\},
$$
where $kl(a,b)$ is the KL divergence between two Bernoulli distributions of respective means $a$ and $b$. Now combining the above inequality to (\ref{eq:low1}) yields a lower bound on a weighted sum of the expected numbers of times each state-control pair is selected. These numbers are directly related to the regret as one can show that:
\begin{align*}
\mathbb{E}_\phi^\pi[R^\pi(T)] =  \sum_{x \in {\cal X}} \sum_{u \notin {\cal O}(x, \phi)}\mathbb{E}_{\phi}^\pi[N_{xu}(T)] \delta^\star(x,u;\phi)  + O(1),
\end{align*}
where ${\cal O}(x, \phi)$ denotes the set of optimal control inputs in state $x$ under $\phi$, and $\delta^\star(x,u;\phi)$ is the sub-optimality gap quantifying the regret obtained by selecting the control $u$ is state $x$. It remains to select the event $E$ or the random variable $Z$ to get a regret lower bound.

To derive problem-specific regret lower bounds, we introduce the notion of {\it uniformly good} algorithms. $\pi$ is uniformly good if for any ergodic MDP $\phi$, any initial state and any constant $\alpha>0$, $\mathbb{E}_\phi^\pi[R^\pi(T)] = o(T^\alpha)$. As it will become clear later, uniformly good algorithms exist. Now select the event $E$ as:\\
$
E = \left[ N_x(T) \ge \rho T,  
 \sum_{u \notin {\cal O}(x, \phi)}  N_{xu}(T) \le  \sqrt{T} \right],
$
for some $\rho>0$ and where $N_x(T)$ is the number of times $x$ is visited up to time $T$. $\rho$ is chosen such that $N_x(T) \ge \rho T$ is very likely under $\pi$ invoking the ergodicity of the MDP. In the change-of-measure argument, $\psi$ is chosen such that $\Pi^\star(\phi)\cap\Pi^\star(\psi)=\emptyset$, where $\Pi^\star(\phi)$ is the set of optimal policies under $\phi$. Now if $\pi$ is uniformly good, $E$ is very likely under $\phi$, and very unlikely under $\psi$. Formally, we can establish that: 
$$
kl(\mathbb{P}_\phi^\pi[E],  \mathbb{P}_\psi^\pi[E])\ge \log(T)+O(1).
$$
Putting the above ingredients together, we obtain:

\begin{theorem}\label{th:low1} (Theorem 1, \cite{ok2018})
Let $\phi$ be an ergodic MDP. For any uniformly good algorithm $\pi$ and for any initial state, 
\begin{align} \label{eq:g-r-lower}
\liminf_{T \to \infty} \frac{\mathbb{E}_\phi^\pi[R^\pi(T)]}{\log T} \ge K(\phi), 
\end{align}
where $K(\phi)$ is the value of the following optimization problem:
\begin{align}
& \underset{\eta \in {\cal F}_\Phi(\phi) }{\textnormal{inf}} 
~~\sum_{x,u} \eta(x, a)  \delta^\star(x, a;\phi), \label{eq:g-r-obj}
\end{align}
where ${\cal F}(\phi)$ is the set of $\eta\ge 0$ satisfying
\begin{equation}
\sum_{(x,u) \in {\cal X} \times {\cal U}}\eta(x, u) KL_{\phi | \psi}(x,u) \ge 1, \ \ \ \forall \psi \in \Delta(\phi)
\end{equation}
and  
$\Delta(\phi) = \{\psi :  \phi \ll \psi,  \Pi^*(\phi) \cap \Pi^*(\psi) = \emptyset \}$.
\end{theorem}

In the above theorem $\psi\ll\phi$ means that the observations under $\psi$ have distributions absolutely continuous w.r.t. those of the observations under $\phi$. By imposing the constraint for all $\psi\in \Delta(\phi)$, we consider all possible confusing MDP $\psi$. It can be shown that the set of  constraints defining ${\cal F}(\phi)$ can be reduced and decoupled: it is sufficient to consider $\psi$ different than $\phi$ in only one suboptimal state-control pair. In that case, the constraints can be written in the following form:  for any $x$ and $u\notin {\cal O}(x,\phi)$, $\eta(x, u) (\delta^\star(x,u;\phi))^2\ge K$ for some absolute constant $K$. As a consequence, the regret lower bound scales as $n_xn_u\log(T)$. The theorem also indicates the optimal exploration rates of the various suboptimal state-control pairs: $\eta(x,u)\log(T)$ represents the expected number of times $(x,u)$ should be observed. Finally, we note that the method used to derive the regret lower bound can also be applied to obtain finite-time regret lower bounds, as in \cite{garivier193}.

For minimax lower bounds, we do not need to restrict the attention to uniformly good algorithms, since for any given algorithm, we are free to pick the MDP for which the algorithm performs the worst. Instead, to identify a ${\cal F}_T^\pi$-measurable $[0,1]$-valued random variable $Z$, such that $kl(\mathbb{E}_\phi^\pi[Z],  \mathbb{E}_\psi^\pi[Z])$ is large, we leverage symmetry arguments. Specifically, The MDP is constructed so as to contain numerous {\it equivalent} states and control inputs, and $Z$ is chosen as the proportion of time a particular state-action pair is selected. Refer to \cite{auer2009near} for the construction of this MDP, and to \cite{garivier193} for more detailed explanations on how to apply change-of-measure arguments to derive minimax bounds.

\begin{theorem}\label{th:low2} (Theorem 5, \cite{auer2009near})
For any algorithm $\pi$, for all integers $n_x, n_u\ge 10$, $D \ge 20\log_A(n_x)$, and $T \ge Dn_xn_u$, there is an MDP with $n_x$ states, $n_u$ control inputs, and diameter $D$ such that for any initial state:
\begin{equation}\label{eq:low2}
\mathbb{E}_\phi^\pi[R^\pi(T)]\ge 0.015 \sqrt{Dn_xn_uT}.
\end{equation}
\end{theorem}

\medskip

\subsubsection{Efficient Algorithms} 

A plethora of learning algorithms have been developed for average-cost MDPs. We can categorize these algorithms based on their design principles. A first class of algorithms aim at matching the asymptotic problem-specific regret lower bound derived above \cite{burnetas1997optimal,graves1997,ok2018}. These algorithms rely on estimating the MDP parameters and in each round (or periodically) they solve the optimization problem (\ref{eq:g-r-obj}) where the true MDP parameters are replaced by their estimators. The solution is then used to guide and minimize the exploration process. This first class of algorithms is discussed further in \S \ref{sec:structure}.

The second class of algorithms includes UCRL, UCRL2 and KL-UCRL\cite{ucrl,auer2009near,filippi}. These algorithms apply the "optimism in front of uncertainty" principle and exhibit finite-time regret ganrantees. They consists in building confidence upper bounds on the parameters of the MDP, and based on these bounds select control inputs. The regret guarantees are anytime, but use worst-case regret as a performance benchmark. For example, UCRL2 with the confidence parameter $\delta$ as an input satisfies:

\begin{theorem} (Theorem 2, \cite{auer2009near})
With probability at least $1-\delta$, for any initial state and all $T\ge 1$, the regret under $\pi= \ $UCRL2 satisfies:
$$
R^\pi(T)\le 34 Dn_x\sqrt{n_uT\log({T\over \delta})}.
$$
\end{theorem}

Anytime regret upper bounds with a logarithmic dependence in the time horizon have been also investigated for UCRL2 and KL-UCRL. For instance, UCRL2 is known to yield a regret in ${\cal O}(D^2n_x^2n_u\log(T/\delta)$ with probability at least $1-3\delta$. 

The last class of algorithms apply a similar Bayesian approach as that used by the celebrated Thompson sampling algorithm for bandit problems. In \cite{agrawal2017}, AJ, a posterior sampling algorithm, is proposed and enjoys the following regret guanrantees:

\begin{theorem} (Theorem 1, \cite{agrawal2017})
With probability at least $1-\delta$, for any initial state, the regret under $\pi= \ $AJ with confidence parameter $\delta$ satisfies: for $T\ge Dn_u \log(T/\delta)^2$,
$$
R^\pi(T)=\tilde{\cal O}(D\sqrt{n_xn_uT}). 
$$
\end{theorem}

\medskip

\subsubsection{Structured MDPs}\label{sec:structure} 

The regret lower bounds derived for tabular MDPs have the major drawback of scaling with the product of the numbers of states and controls, $n_xn_u$. Hence, with large state and control spaces, it is essential to identify and exploit any possible structure existing in the system dynamics and cost function so as to minimize exploration phases and in turn reduce regret to reasonable values. Modern RL algorithms actually implicitly impose some structural properties either in the model parameters (transition probabilities and cost function, see e.g. \cite{ortner2012online}) or directly in the $Q$-function (for discounted RL problems, see e.g. \cite{mnih015}. Despite their successes, we do not have any regret guarantees for these recent algorithms. Recent efforts to develop algorithms with guarantees for structured MDPs include \cite{ortner2012online,lakshmanan2015improved,eluder,ok2018}. \cite{ok2018} is the first paper extending the analysis of  \cite{burnetas1997optimal} to the case of structured MDPs. The authors derive a problem-specific regret lower bound, and show that the latter can be obtained by just modifying in Theorem \ref{th:low1} the definition of the set of confusing MDPs $\Delta(\phi)$. Specifically, if $\Phi$ denotes a set of structured MDPs ($\Phi$ encodes the structure), then $\Delta(\phi) = \{\psi \in \Phi:  \phi \ll \psi,  \Pi^*(\phi) \cap \Pi^*(\psi) = \emptyset \}$. The minimal expected regret scales as $K_\Phi(\phi)\log(T)$, where $K_\Phi(\phi)$ is the value of the modified optimization problem (\ref{eq:g-r-obj}). In \cite{ok2018}, DEL, an algorithm extending that proposed in \cite{burnetas1997optimal} to the case of structured MDPs, is shown to optimally exploit the structure:

\begin{theorem} (Theorem 4, \cite{ok2018})
For any $\phi\in \Phi$, the regret under $\pi= \ $DEL satisfies:
$$
\lim\sup_{T\to\infty} {\mathbb{E}_\phi^\pi[R^\pi(T)]\over \log(T)} \le K_\Phi(\phi).
$$
\end{theorem}
 
The semi-infinite LP characterizing the regret lower bound can be simplified for some particular structures. For example in the case where $p$ and $q$ smoothly vary over states and controls (Lipschitz continuous), it can be shown that the regret lower bound does not scale with $n_x$ and $n_u$. The simplified LP can then be used as in \cite{burnetas1997optimal} to devise an asymptotically optimal algorithm.

\subsection{Discounted MDPs}  

Most research efforts, from early work \cite{watkins1992q} to more recent Deep RL \cite{mnih015}, towards the design of efficient algorithms for discounted MDPs have focussed on model-free approaches, where one directly learns the value or the Q-value function of the MDP. Such an approach leads to simple algorithms that are potentially more robust than model-based algorithms (since they do not rely on modelling assumptions). The performance analysis of these algorithms has been initially mainly centered around the question of their convergence; for example, the analysis of $Q$-learning algorithm with function approximation \cite{maei2009convergent} often calls for new convergence results of stochastic approximation schemes. Researchers have then strived to investigate and optimize their convergence rates. There is no consensus on the metric one should use to characterize the speed of convergence; e.g. the recent Zap Q-learning algorithm \cite{Devraj2017zap} minimizes the asymptotic error covariance matrix, while most other analayses focus on the minimax sample complexity. Note that the notion of regret in discounted settings is hard to define and has hence not been studied. Also observe that problem-specific metrics have not been investigated yet, and it hence seems perilous to draw definitive conclusions from existing theoretical results for learning discounted MDPs. 

\subsubsection{Sample complexity lower bound} For discounted MDPs, the sample complexity is defined as the number of samples one need to gather so as to learn an $\epsilon$-optimal policy with probability at least $1-\delta$. Minimax sample complexity lower bounds are known for both the generative and online sampling models:

\begin{theorem} (Theorem 1, \cite{azar2012sample} and Theorem 11, \cite{lattimore2014near}) In both the generative and online sampling models, for $\epsilon$ and $\delta$ small enough, there exists an MDP such that any learning algorithm must have a sample complexity in $\Omega({n_xn_u\over \epsilon^2(1-\gamma)^3}\log({n_x\over \delta}))$ (where $\gamma$ denotes the discount factor).
\end{theorem}

\subsubsection{The price of model-free approaches}

Some model-based algorithms are known to match the minimax sample complexity lower bound. In the online sampling setting, the authors of \cite{lattimore2014near} presents UCRL($\gamma$), an extension of UCRL for discounted costs, and establish a minimax sample complexity upper bound matching the above lower bound. UCRL($\gamma$) consists in deriving upper confidence bounds for the MDP parameters, and in selecting action optimistically (this can lead to important computational issues). In the generative sampling model, algorithms mixing model-based and model-free approaches have been shown to be minimax-sample optimal. This is the case of QVI (Q-value Iteration) initially proposed in \cite{kearns2002near} and analyzed in \cite{azar2012sample}. QVI estimates the MDP, and from this estimator, applies a classical value iteration method to approximate the $Q$-function and hence the optimal policy. QVI can be also made computationally efficient \cite{sidford2018}.  

As for now, there is no pure model-free algorithm achieving the minimax sample complexity limit. Speedy Q-learning \cite{azar2011speedy} has a minimax sample complexity in $\tilde{\cal O}({n_xn_u\over \epsilon^2(1-\gamma)^4})$ (this is for now the best one can provably do using model-free approaches). However, there is hope to find model-free minimax-sample optimal algorithms. In fact, recently, $Q$-learning with exploration driven by simple upper confidence bounds on the $Q$-values (rather than on the MDP parameters as in UCRL) has been shown to be minimax regret optimal for episodic reinforcement learning tasks \cite{jin2018}. It is likely that model-free algorithms can be made minimax optimal. If this is verified, this would further advocate the use of model-specific rather than minimax performance metrics.


\section{Model-Based Methods for LQR}
\label{sec:model_based_lqr}
We combine the techniques described in \cite{sysID} with robust and optimal control to derive finite-time guarantees for the optimal LQR control of an unknown system.  We partition our study according to three initial uncertainty regimes: (i) completely unknown $(A,B)$, (ii) moderate error bounds under which CE control may fail, and (iii) small error bounds under which CE control is stabilizing.

\subsection{PAC Bounds for Unknown $(A,B)$}
Here we assume that the system is completely unknown, and consider the problem of identifying system estimates $(\Ahat,\Bhat)$, bounding the corresponding parameter uncertainties $\epsilon_A = \norm{\Ahat - A}_2$ and $\epsilon_B = \norm{\Bhat - B}_2$, and using these system estimates and uncertainty bounds to compute a controller with provable performance bounds.  In what follows, unless otherwise specified, all results are taken from \cite{learning-lqr}.  

The system identification and uncertainty quantification steps are covered in Theorem IV.3 of \cite{sysID}, which we summarize here for the convenience of the reader.

Consider a linear dynamical system described by
\begin{equation}
x_{t+1} = Ax_t + Bu_t + w_t,
\label{eq:model-ls}
\end{equation}
where $x_t,w_t \in \R^{n_x}$, $u_t \in \R^{n_u}$, and $w_t \iid \Normal(0,\sigma_w^2 I)$.  To identify the matrices $(A,B)$,
we inject excitatory Gaussian noise via $u_t \iid \Normal (0,\sigma_u^2 I)$.  We run $N$ experiments over a horizon of $T+1$ time-steps, and then solve for our estimates $(\Ahat,\Bhat)$ via the least-squares problem:
\begin{equation}\begin{array}{rcl}
\begin{bmatrix} \Ahat & \Bhat \end{bmatrix}^\top &=& \arg\displaystyle\min_{(A,B)} \sum_{i=1}^N \norm{x_{T+1}^{(i)} - Ax_{T}^{(i)}- Bu_T^{(i)}}_2^2.
\end{array}
\label{eq:matrix-ls}
\end{equation}
Notice that we only use the last time-steps of each trajectory: we do so for analytic simplicity, and return to single trajectory estimators that use all data later in the section.  We then have the following guarantees.

\begin{theorem}
Consider the least-squares estimator defined by \eqref{eq:matrix-ls}.  Fix a failure probability $\delta \in (0,1]$, and assume that $N \geq 24(n_x + n_u)\log(54/\delta)$.  Then, it holds with probability at least $1-\delta$, that
\begin{align}
\twonorm{\Ahat - A} \leq {8\sigma_w}{{\lambda^{-1/2}_{\min}\left(\Sigma_x\right)}}\sqrt{\frac{(2n_x+n_u)\log(54/\delta)}{N}},
\label{eq:Aerror}
\end{align}
and
\begin{equation}
\twonorm{\Bhat - B} \leq \frac{8\sigma_w}{\sigma_u}\sqrt{\frac{(2n_x + n_u)\log(54/\delta)}{N}},
\label{eq:Berror}
\end{equation}
where $\lmin{\Sigma_x}$ is the minimum eigenvalue of the finite time controllability Gramian $\Sigma_x := \sigma_u^2 \sum_{t=0}^T A^tBB^\top (A^\top)^t + \sigma_w^2 \sum_{t=0}^T A^t (A^\top)^t.$
\label{thm:iid-matrix-error}
\end{theorem}

We now condition on the high-probability guarantee Theorem \ref{thm:iid-matrix-error}, and assume that we have system estimates $(\Ahat,\Bhat)$ and corresponding uncertainty bounds $(\epsilon_A,\epsilon_B)$, allowing us to focus on the controller synthesis step.  Our goal is to compute a controller that is robustly stabilizing for any admissible realization of the system parameters, and for which we can bound performance degradation as a function of the uncertainty sizes $(\epsilon_A, \epsilon_B)$.

In order to meet these goals, we use the \emph{System Level Synthesis} \cite{SysLevelSyn1,anderson2019} (SLS) nominal and robust \cite{matni2017scalable} parameterizations of stabilizing controllers.  The SLS framework focuses on the \emph{system responses} of a closed-loop system.  Consider a LTI causal controller $\tf K$, and let $\tf u = \tf K \tf x$. Then the closed-loop transfer matrices from the process noise $\tf w$ to the state $\tf x$ and control action $\tf u$ satisfy

\begin{equation}
\begin{bmatrix} \tf x \\ \tf u \end{bmatrix} = \begin{bmatrix} (zI - A-B\tf K)^{-1} \\ \tf K (zI-A-B \tf K)^{-1} \end{bmatrix} \tf w.
\label{eq:response}
\end{equation}
We then have the following theorem parameterizing the set of stable closed-loop transfer matrices, as described in equation \eqref{eq:response}, that are achievable by a stabilizing controller $\tf K$.
\begin{theorem}[State-Feedback Parameterization~\cite{SysLevelSyn1}]
The following are true:
\begin{itemize}
\item The affine subspace defined by
\begin{equation}
\begin{bmatrix} zI - A & - B \end{bmatrix} \begin{bmatrix} \tf \Phi_x \\ \tf \Phi_u \end{bmatrix} = I, \ \tf \Phi_x, \tf \Phi_u \in \frac{1}{z}\mathcal{RH}_\infty
\label{eq:achievable}
\end{equation}
parameterizes all system responses \eqref{eq:response} from $\tf w$ to $(\tf x, \tf u)$, achievable by an internally stabilizing state-feedback controller $\tf K$.
\item For any transfer matrices $\{\tf \Phi_x, \tf \Phi_u\}$ satisfying \eqref{eq:achievable}, the controller $\tf K = \tf \Phi_u \tf \Phi_x^{-1}$ is internally stabilizing and achieves the desired system response \eqref{eq:response}.
\end{itemize}
\label{thm:param}
\end{theorem}

We will also make use of the following robust variant of Theorem \ref{thm:param}.
\begin{theorem}[Robust Stability~\cite{matni2017scalable}]
Let $\tf\Phi_x$ and $\tf \Phi_u$ be two transfer matrices in $\frac{1}{z}\mathcal{RH}_\infty$ such that 
\begin{equation}
\begin{bmatrix} zI - A & - B \end{bmatrix} \begin{bmatrix} \tf \Phi_x \\ \tf \Phi_u \end{bmatrix} = I + \tf \Delta.
\end{equation}
Then the controller $\tf K = \tf \Phi_u \tf \Phi_x^{-1}$ stabilizes the system described by $(A,B)$ if and only if $(I+\tf \Delta)^{-1} \in \mathcal{RH}_\infty$.  Furthermore, the resulting system response is given by
\begin{equation}
\begin{bmatrix} \tf x \\ \tf u \end{bmatrix} = \begin{bmatrix} \tf \Phi_x \\ \tf \Phi_u \end{bmatrix}(I+\tf \Delta)^{-1} \tf w.
\end{equation}
\label{thm:robust}
\end{theorem}
\begin{coro}
Under the assumptions of Theorem \ref{thm:robust}, if $\|\tf \Delta \| <1$ for any induced norm $\|\cdot \|$, then the controller $\tf K = \tf \Phi_u \tf \Phi_x^{-1}$ stabilizes the system described by $(A,B)$.
\label{coro:sufficient}
\end{coro}

We now return to the problem setting where the estimates $(\Ah, \Bh)$ of a true system $(A,B)$ satisfy $\|\Delta_A\|_2\leq\epsilon_A,~~\|\Delta_B\|_2\leq\epsilon_B$, for  $\Delta_A := \Ah-A$ and $\Delta_B := \Bh-B$.  We first formuate the LQR problem in terms of the system responses $\{\tf\Phi_x,\tf\Phi_u\}$. It follows from Theorem \ref{thm:param} and the standard equivalence between infinite horizon LQR and $\htwo$ optimal control that, for a disturbance process distributed as $w_t \iid \mathcal{N}(0,\sigma_w^2 I)$, the standard LQR problem \eqref{eq:LQR-T} can be equivalently written as
\begin{equation}
\min_{\tf\Phi_x, \tf\Phi_u} \sigma_w^2 \left\|\begin{bmatrix} Q^\frac{1}{2} & 0 \\ 0 & R^\frac{1}{2}\end{bmatrix}\begin{bmatrix} \tf \Phi_x \\ \tf \Phi_u \end{bmatrix}\right\|_{\htwo}^2 \text{ s.t. equation \eqref{eq:achievable}}.
\label{eq:lqr2}
\end{equation}
Going forward, we drop the $\sigma_w^2$ multiplier in the objective function as it does not affect the guarantees that we compute.

We begin with a simple sufficient condition under which any controller $\tf K$ that stabilizes $(\Ah,\Bh)$ also stabilizes the true system $(A,B)$.  For a matrix $M$, we let $\Res{M}$ denote the resolvent, i.e., $\Res{M} := (zI - M)^{-1}\,.$

\begin{lemma}\label{lemma:robust-sls}
Let the controller $\tf K$ stabilize $(\Ah, \Bh)$ and $(\tf\Phi_x,\tf\Phi_u)$ be its corresponding system response \eqref{eq:response} on system $(\Ah,\Bh)$.  Then if $\tf K$ stabilizes $(A,B)$, it achieves the following LQR cost $J(A,B,\tf K)$ defined as
\begin{equation}
 {\scriptsize \left\|\begin{bmatrix} Q^\frac{1}{2} & 0 \\ 0 & R^\frac{1}{2}\end{bmatrix}\begin{bmatrix} \tf\Phi_x \\  \tf\Phi_u \end{bmatrix}\left(I+\begin{bmatrix}\Delta_A& \Delta_B\end{bmatrix}\begin{bmatrix} \tf\Phi_x \\  \tf\Phi_u \end{bmatrix}\right)^{-1}\right\|_{\htwo}}\:.
\end{equation}

Furthermore, letting
\begin{equation}\label{eq:deltahat}
\tf{\hat\Delta} := \begin{bmatrix}\Delta_A& \Delta_B\end{bmatrix}\begin{bmatrix} \tf\Phi_x \\  \tf\Phi_u \end{bmatrix} = (\Delta_A + \Delta_B \tf K)\Res{\Ah+\Bh\tf K} \:.
\end{equation}
controller $\tf K$ stabilizes $(A,B)$ if $\hinfnorm{\tf{\hat{\Delta}}} <1$.
\label{lem:sufficient}
\end{lemma}
We can therefore pose a robust LQR problem as
\begin{align}\label{eq:robustLQR}
\begin{split}
 &\min_{\tf\Phi_x, \tf \Phi_u} \sup\limits_{\substack{\|\Delta_A\|_2\leq \epsilon_A \\ \|\Delta_B\|_2\leq \epsilon_B}}  J(A,B,\tf K)\\
& \text{s.t.} \begin{bmatrix}zI-\Ah&-\Bh\end{bmatrix}\begin{bmatrix} \tf\Phi_x \\  \tf\Phi_u \end{bmatrix} = I,~~\tf\Phi_x, \tf \Phi_u  \in\frac{1}{z}\mathcal{RH}_\infty \:.
 \end{split}
\end{align}
The resulting robust control problem is one subject to real-parametric
uncertainty, a class of problems known to be computationally
intractable~\cite{braatz94}.  To circumvent this issue, we instead find an upper-bound to the cost $J(A,B,\tf K)$ that is independent of the uncertainties $\Delta_A$ and $\Delta_B$. First, note that if $\hinfnorm{\Dh} < 1$, we can write
\begin{align} \label{eq:upperbnd1}
J(A,B,\tf K) 
&\leq \frac{1}{1-\hinfnorm{\Dh}}J(\Ah,\Bh,\tf K).
\end{align}

This upper bound separates nominal performance, as captured by $J(\Ah,\Bh,\tf K)$, from the effects of the model uncertainty, as captured by $(1-\hinfnorm{\Dh})^{-1}$.  It therefore remains to compute a tractable bound for $\hinfnorm{\Dh}$.
\begin{proposition}[Proposition 3.5, \cite{learning-lqr}]
	For any $\alpha \in (0,1)$ and $\tf{\hat{\Delta}}$ as defined in \eqref{eq:deltahat}, we have
\begin{equation}\label{eq:ben-tri-bound}
	\|  \tf{\hat{\Delta}} \|_{\hinf}
	 \leq \left\|\begin{bmatrix} \tfrac{\epsilon_A}{\sqrt{\alpha}} \tf \Phi_x \\ \tfrac{\epsilon_B}{\sqrt{1-\alpha}} \tf\Phi_u \end{bmatrix} \right\|_{\hinf} =\colon H_\alpha(\tf\Phi_x,\tf\Phi_u) \:.
\end{equation}
\label{prop:bound}
\end{proposition}
Applying Proposition \ref{prop:bound} in conjunction with the bound \eqref{eq:upperbnd1}, we arrive at the following upper bound to the cost function of the robust LQR problem \eqref{eq:robustLQR}, which is independent of the perturbations $(\Delta_A,\Delta_B)$:
\begin{align}\label{eq:upperbnd}
 \sup\limits_{\|\Delta_A\|_2\leq \epsilon_A, \, \|\Delta_B\|_2\leq \epsilon_B}  J(A,B,\tf K)  &\leq \frac{J(\Ah,\Bh,\tf K)}{1 - H_\alpha(\tf\Phi_x,\tf\Phi_u)}\:.
\end{align}
The upper bound is only valid when $H_\alpha(\tf \Phi_x, \tf \Phi_u) < 1$, which guarantees the stability of the closed-loop system.  Note that \eqref{eq:upperbnd} can be used to upper bound the performance achieved by any robustly stabilizing controller.

We can then pose the robust LQR synthesis problem as the following quasi-convex optimization problem, which can be solved by gridding over $\gamma\in[0,1)$:
\begin{align}\label{eq:robustLQRbnd}
\begin{split}
& \minimize_{\gamma\in[0,1)}\frac{1}{1 - \gamma}\min_{\tf\Phi_x, \tf \Phi_u\in\frac{1}{z}\mathcal{RH}_\infty} \left\|\begin{bmatrix} Q^\frac{1}{2} & 0 \\ 0 & R^\frac{1}{2}\end{bmatrix}\begin{bmatrix} \tf\Phi_x \\  \tf\Phi_u \end{bmatrix}\right\|_{\htwo}\\&
 \text{s.t.} \begin{bmatrix}zI-\Ah&-\Bh\end{bmatrix}\begin{bmatrix} \tf\Phi_x \\  \tf\Phi_u \end{bmatrix} = I,~~\left\|\begin{bmatrix} \tfrac{\epsilon_A}{\sqrt{\alpha}} \tf \Phi_x \\ \tfrac{\epsilon_B}{\sqrt{1-\alpha}} \tf\Phi_u \end{bmatrix} \right\|_{\hinf}\leq \gamma\\
 \end{split}
\end{align}
As we constrain $\gamma \in [0,1)$, any feasible solution $(\tf \Phi_x, \tf \Phi_u)$ to optimization problem \eqref{eq:robustLQRbnd} generates a controller $\tf K = \tf \Phi_u \tf \Phi_x^{-1}$ that stabilizes the true system $(A,B)$.  

\begin{remark}
Optimization problem \eqref{eq:robustLQRbnd} is infinite-dimensional.  However, one can solve a finite-dimensional approximation of the problem over a horizon $T = \Omega(\log(1/(\epsilon_A + \epsilon_B))$ (see Theorem 5.1, \cite{learning-lqr}) such that the sub-optimality bounds we prove below still hold up to universal constants. 
\end{remark}  

We then have the following theorem bounding the sub-optimality of the proposed robust LQR controller.


\begin{theorem}
\label{thm:lqr_cost}
Let $J_\star$ denote the minimal LQR cost achievable by any controller for the dynamical system with transition matrices $(A,B)$, and let $\trueK$ denote the optimal contoller. Let $(\Ah,\Bh)$ be estimates of the transition matrices such that $\ltwonorm{\Delta_A} \leq \epsilon_A$, $\ltwonorm{\Delta_B} \leq \epsilon_B$. Then, if $\tf K$ is synthesized via \eqref{eq:robustLQRbnd} with $\alpha = 1/2$, the relative error in the LQR cost is
\begin{align} \label{eq:lqr_bound}
\frac{J(\trueA, \trueB, \tf K) - J_\star }{J_\star} \leq 5 (\epsilon_A + \epsilon_B\ltwonorm{\trueK})\hinfnorm{\Res{A+B\trueK}} \:,
\end{align}
as long as $(\epsilon_A + \epsilon_B\ltwonorm{\trueK})\|\Res{A+B\trueK}\|_{\hinf}\leq 1/5$.
\end{theorem}

The crux of the proof of Theorem \ref{thm:lqr_cost} is to show that for sufficiently small $(\epsilon_A,\epsilon_B)$, the optimal LQR controller $\trueK$ is robustly stabilizing.  We then exploit optimality to upper bound the cost $(1-\gamma)^{-1} J(\Ahat,\Bhat,\tf K)$ achieved by our controller with that achieved by the optimal LQR controller $(1-\gamma_{LQR})^{-1}J(\Ahat,\Bhat,\trueK)$, from which the result follows almost immediately by repeating the argument with estimated and true systems reversed.  Combining Theorems \ref{thm:lqr_cost} and \ref{thm:iid-matrix-error}, we see that  ${J(A,B,\tf K) - J_\star\leq O(\epsilon_A+\epsilon_B) \leq O(\sqrt{(n+p)\log(1/\delta)/N})}$.  This in turn shows that LQR optimal control of an unknown system is $(\epsilon,\delta)$-episodic PAC learnable, where here we interpret each system identification experiment as an episode.

\begin{example}
 \begin{figure}[!t]
 \centering
\includegraphics[width=.475\columnwidth]{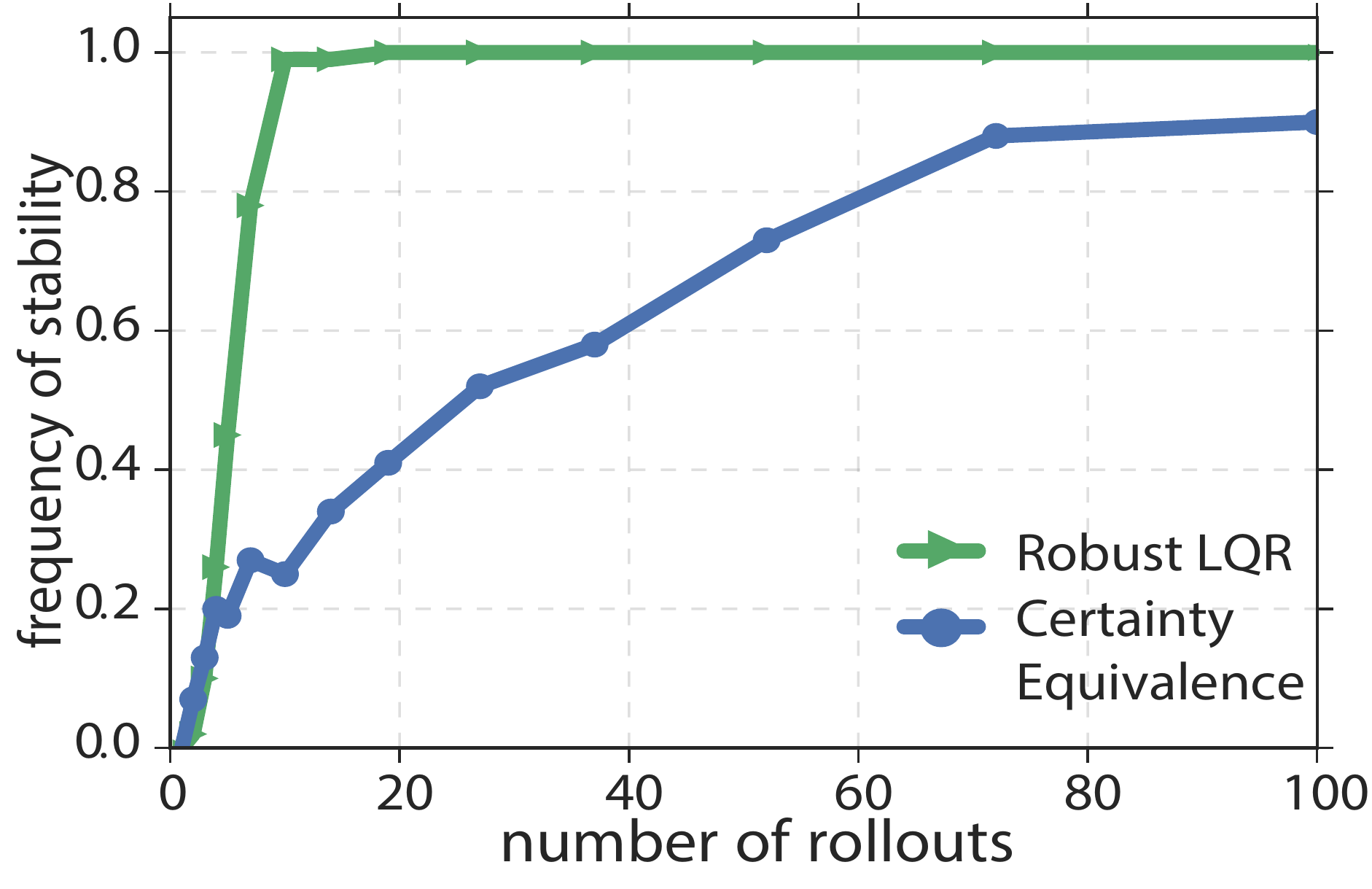}~
\includegraphics[width=.475\columnwidth]{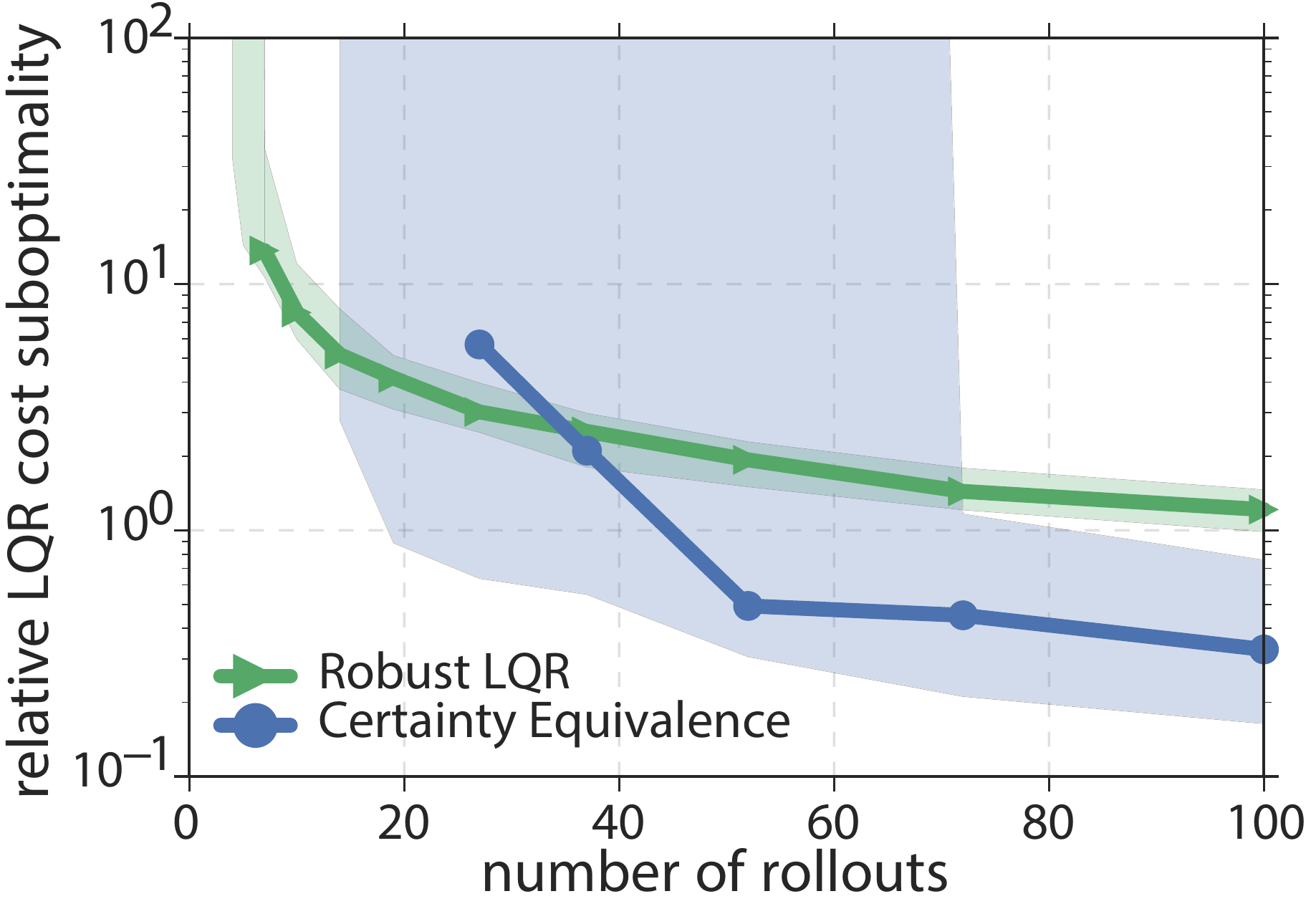}
\caption{\small Left: The percentage of stabilizing controllers synthesized using certainty equivalence and robust LQR over 100 independent trials.  Model estimates and uncertainty bounds are computed using independent rollouts with horizon $T=6$.  Right: Corresponding sub-optimality bounds.  CE controllers outperform robust LQR controllers when they are stabilizing.}
\label{fig:stability}
\end{figure}
Consider an LQR problem specified by
\begin{align} \label{eq:exampledynamics}
A  = \begin{bmatrix} 1.01 & 0.01 & 0\\
0.01 & 1.01 & 0.01\\
0 & 0.01 & 1.01\end{bmatrix},\, B = I,\, Q = 10^{-3} I,\, R =I \:.
\end{align}
In the left plot of Figure \ref{fig:stability}, we show the percentage of stabilizing controllers synthesized using certainty equivalence and the proposed robust LQR controllers over 100 independent trials..  Notice that even after collecting data from 100 trajectories, the CE controller yields unstable behavior in approximately 10\% of cases.  Given that the state of the underlying dynamical system is only 3-dimensional, one might consider 100 data-points to be a reasonable approximation of an ``asymptotic'' amount of data, highlighting the need for a more refined analysis of the effects of finite data on stability and performance.  Contrast this with the behavior achieved by the robust LQR synthesis method, which explicitly accounts for system uncertainty: after a small number of trials, there is a sharp transition to 100\% stability across trials.  Further, feasibility of the synthesis problem \eqref{eq:robustLQRbnd} provides a certificate of stability and performance, conditioned on the uncertainty bounds being correct.  However, robustness does come at a price: as shown in the right plot of Figure \ref{fig:stability}, the CE controller outperforms the robust LQR controller when it is stabilizing.
\end{example}

\subsection{Regret Bounds under Moderate Uncertainty}
We have just described an offline procedure for learning a coarse estimate of system dynamics and computing a robustly stabilizing controller.  We now consider the task of adaptively refining this model and controller.  For this problem, we seek high probability bounds on the regret $R(T)$, defined as
\begin{equation}
R(T) := \sum_{t=0}^T (x_t^\top Q x_t + u_t^\top R u_t) - TJ_\star,
\label{eq:lqr_regret}
\end{equation}  
for $J_\star$ defined as in the previous section.

Consider the single trajectory least-squares estimator:
\begin{align}
	(\Ah, \Bh) \in \arg\min_{A,B} \frac{1}{2} \sum_{t=0}^{T-1} \norm{x_{t+1} - A x_t - B u_t}^2_2 \:, \label{eq:single_trajectory_ls}
\end{align}
solved with data generated from system \eqref{eq:model-ls} driven by input $u_t \iid \Normal(0,\sigma_u^2 I_{n_u})$.
The following result from Simchowitz et al.~\cite{simchowitz2018learning} gives us a
high probability bound on the error of the estimator \eqref{eq:single_trajectory_ls}.
\begin{theorem}[\cite{simchowitz2018learning}]
\label{thm:single_trajectory_ls}
Suppose that $A$ is stable
and that the trajectory length $T$ satisfies:
\begin{align*}
	T \geq \Omega(n_x + n_u + n_x \log(\frac{1}{\delta} (1 + \sigma_u^2 \norm{B}^2 / \sigma_w^2))) \:.
\end{align*}
With probability at least $1-\delta$, the quantity $\max\{ \norm{\Ah - A}, \norm{\Bh - B} \}$ is bounded above by 
\begin{align*}
	O\left( \sigma_w \sqrt{  \frac{ n_x+n_u + n_x \log\left(\frac{1}{\delta} (1 + \frac{\sigma_u^2}{\sigma_w^2} \norm{B}^2) \right) }{ T \min\{  \sigma_w^2, \sigma_u^2  \}}      }    \right) \:.
\end{align*}
Here, the $O(\cdot)$ hides specific properties of the controllability gramian.
\end{theorem}

Suppose we are provided with an initial stabilizing controller $\tf K$: how should we balance controlling the system (exploitation) with exciting it for system identification purposes (exploration)?  We propose studying the simple exploration scheme $\tf u = \tf K \tf x + \tf \eta,$
where $\eta_t \iid \Normal(0,\sigma_\eta^2I_{n_u})$.  Theorem \ref{thm:single_trajectory_ls} tells us that if we collect data for $T$ time-steps and compute estimates $(\Ahat(T),\Bhat(T))$ using ordinary least squares, then with high probability we have that 
\begin{equation}
\epsilon:= \max(\norm{\Ahat(T) - A}_2, \norm{\Bhat(T) - B}_2) \leq \tilde O\left(\frac{1}{\sigma_\eta T^{1/2}}\right) 
\label{eq:adaptive_eps}
\end{equation}  
when $\sigma_\eta \ll \sigma_w$.
Furthermore, we saw in Theorem \ref{thm:lqr_cost} that for a model error size bounded by $\epsilon$, that the sub-optimality incurred by a robust LQR controller $\tf K$ synthesized using problem \eqref{eq:robustLQRbnd} satisfies $\hat J - J_\star \leq \tilde O(\epsilon)$, where here we use $\hat J$ to denote the cost achieved by the controller $\tf K$.  Letting $\hat J_T$ denote the $T$ horizon cost $\E[\sum_{t=1}^T x_t^\top Q x_t + u_t^\top R u_t]$ achieved by the robust LQR controller, we can show using a similar argument that $\hat J_T - TJ_\star \leq \tilde O(T \epsilon)$.  However, we must also consider the performance degradation incurred by injecting the exploratory signal $\tf \eta$. As this signal is independent of all others, it is easy to see that it incurs an additional cost of $\tilde O(T\sigma_\eta^2)$.  Combining these arguments, we conclude that
\begin{equation}
\hat J_T - T J_\star \leq \tilde O\left(\frac{T^{1/2}}{\sigma_\eta}\right) + \tilde O(T\sigma_\eta^2),
\label{eq:explore_exploit}
\end{equation}
where the first term comes from combining the bound $\hat J_T - TJ_\star \leq \tilde O(T \epsilon)$ with \eqref{eq:adaptive_eps}.  Setting $\sigma_\eta^2 = C_\eta T^{-1/3}$ optimizes the right-hand side of the bounds, leading to $\hat J_T - T J_\star \leq \tilde O(T^{2/3})$.

This reasoning was used in \cite{dean2018regret} to develop an algorithm that achieves $\tilde O(T^{2/3})$ regret, with high probability, as captured by the regret measure \eqref{eq:lqr_regret}.  We informally summarize the algorithm and main results below before commenting on the strengths and weaknesses of the method.

\begin{center}
\begin{algorithm}[h!]\scriptsize
   \caption{Robust Adaptive LQR (Informal)}
\begin{algorithmic}[1]
\STATE \textbf{Input:} initial stabilizing controller $\tf K^0$, failure probability $\delta \in (0,1]$, base epoch length $C_T$, base exploration variance $C_\eta$

\FOR{$i=0,1,2,\dots$}
\STATE Set $T_i \leftarrow C_T2^i$, $\sigma_{\eta,i}^2 \leftarrow C_\eta T_i^{-1/3}$
\STATE Collect data $\{x_t^i, u_t^i\}_{t=0}^{T_i} \leftarrow$ evolve system for $T_i$ stps with $\tf u = \tf K^i \tf x + \tf \eta_i$, with $\eta_{i,t}\iid\Normal(0,\sigma_{\eta,i}^2 I_{n_u})$
\STATE $(\hat{A}_i, \hat{B}_i, \epsilon_i) \leftarrow$ solve OLS problem using collected data and estimate uncertainty $\epsilon_i$
\STATE $\tf K^i \leftarrow \mathrm{RobustLQR}(\hat{A}_i,\hat{B}_i,\epsilon_i)$
\ENDFOR
\end{algorithmic}
   \label{alg:robust_adaptive}
 \end{algorithm}
\end{center}
\begin{theorem}[Informal, Theorems 3.2 \& 3.3, \cite{dean2018regret}]
With the system driven by Algorithm \ref{alg:robust_adaptive}, we have with probability at least $1-\delta$ that the estimates at time $T$ satisfy ${\max(\norm{\Ahat-A}_2,\norm{\Bhat-B}_2) \leq \tilde O((n_x+n_u)^{\tfrac{1}{2}}T^{-\tfrac{1}{3}})}$, and that the regret \eqref{eq:lqr_regret} satisfies $R(T) \leq \tilde O((n_x + n_u)T^{2/3})$.
\label{thm:informal_regret}
\end{theorem}

Theorem \ref{thm:informal_regret} tells us that Algorithm \ref{alg:robust_adaptive} simultaneously leads to consistent estimates of the system matrices $(A,B)$ and near optimal performance, while providing high-probability guarantees on the stability and performance of the system over time.  However, as shown in Figure \ref{fig:stability}, adding robustness to model uncertainty incurs a loss of performance.  One might then ask how sub-optimal the $\tilde O(T^{2/3})$ regret bound is.  Indeed, from the linear bandits literature \cite{rusmevichientong2010linearly}, we know that it can be no lower than $\tilde O(T^{1/2})$.  In the next subsection, we summarize the results of \cite{mania2019certainty}, which state that for sufficiently small uncertainty $\epsilon$, CE control is nearly optimal.

\subsection{Regret Bounds under Small Uncertainty}
\label{sec:nominal}

Suppose that we can show that a nominal controller $K=K(\Ahat,\Bhat)$ computed using model estimates $(\Ahat,\Bhat)$ satisfying error bound \eqref{eq:adaptive_eps}, i.e., $\max(\norm{\Ahat-A}_2,\norm{\Bhat-B}_2) \leq \epsilon$,  satisfies $\norm{K - \trueK} \leq O(\epsilon)$.  Then, if we were able to take a 2nd order Taylor series expansion of the LQR cost $J(K)$ around the optimal controller $\trueK$, we would observe that:
\begin{align*}
J(K) &= J(\trueK)+ \langle \nabla J(\trueK),K-\trueK\rangle \\
&\quad + \frac{1}{2}\langle K-\trueK, \nabla^2 J\left((1-\gamma)\trueK + \gamma K\right)\cdot\left(K-\trueK\right)\rangle \\
&= J(\trueK) + O(\epsilon^2),
\end{align*}
where the first equality holds for some $\gamma \in (0,1)$ by the mean value form of the Taylor series expansion, and the second equality holds by recognizing that $\nabla J(\trueK) = 0$ as it must be a stationary point of the cost functional $J$.  This intuitive argument is formalized in \cite{mania2019certainty}, wherein they explicitly quantify a bound on the error $\epsilon$ such that this approximation is valid.

\begin{theorem}[Informal, Theorem 2, \cite{mania2019certainty}]
Let $\epsilon>0$ be such that $\norm{\Ahat-A}_2,\norm{\Bhat-B}_2 \leq \epsilon$, and assume that $Q,R \succ 0$.  Then the cost $\hat J$ achieved by applying control input $u_t = K(\Ahat,\Bhat)x_t$ satisfies
\begin{equation}
\hat J - J_\star \leq O\left( n_u n_x^5\epsilon^2\right)
\label{eq:ce_bounds}
\end{equation}
so long as $\epsilon$ is sufficiently small.
\label{thm:ce_bound}
\end{theorem}

In the interest of space, we do not expand all of the problem dependent constants in bound \eqref{eq:ce_bounds}, nor do we explicitly describe the required bounds on $\epsilon$; these are however available in \cite{mania2019certainty}.  It is however worth noting that the bound \eqref{eq:ce_bounds} is applicable to a much smaller size of model uncertainty than the robust bound \eqref{eq:lqr_bound} provided in Theorem \ref{thm:lqr_cost}.
Within this local neighborhood, Theorem \ref{thm:ce_bound} implies that CE control leads to performance satisfying $\hat J - J_\star \leq O(\epsilon^2)$, whereas the performance of the robust LQR controller only satisfies $\hat J - J_\star \leq O(\epsilon)$.  Further, integrating this bound into the exploration/exploitation tradeoff \eqref{eq:explore_exploit}, we see that the right hand side is now minimized by setting $\sigma_\eta^2 = C_\eta T^{-1/2}$, leading to the desired $\hat J_T - T J_\star \leq \tilde O(T^{1/2})$.

\begin{example}
We consider the same dynamics \eqref{eq:exampledynamics}, but set $Q=10I$.  In Figure \ref{fig:regrets}, from \cite{dean2018regret}, we show a comparison of different adaptive methods on 500 experiments.  The median and 90th percentiles are shown for the optimal controller, the certainty equivalence controller, the robust LQR controller, and heuristic implementations of the Thompson Sampling (TS) method proposed in \cite{ouyang17} and the Optimism in the Face of Uncertainty (OFU) based method proposed in \cite{abbasi2011regret}.  All methods achieve similar performance, and we note that the latter two methods assume that CE like control is guaranteed to be stabilizing, and further, that the OFU method requires solving a non-convex optimization problem as a subroutine.
\begin{figure}[!t]
\centering
\includegraphics[width=.5\columnwidth]{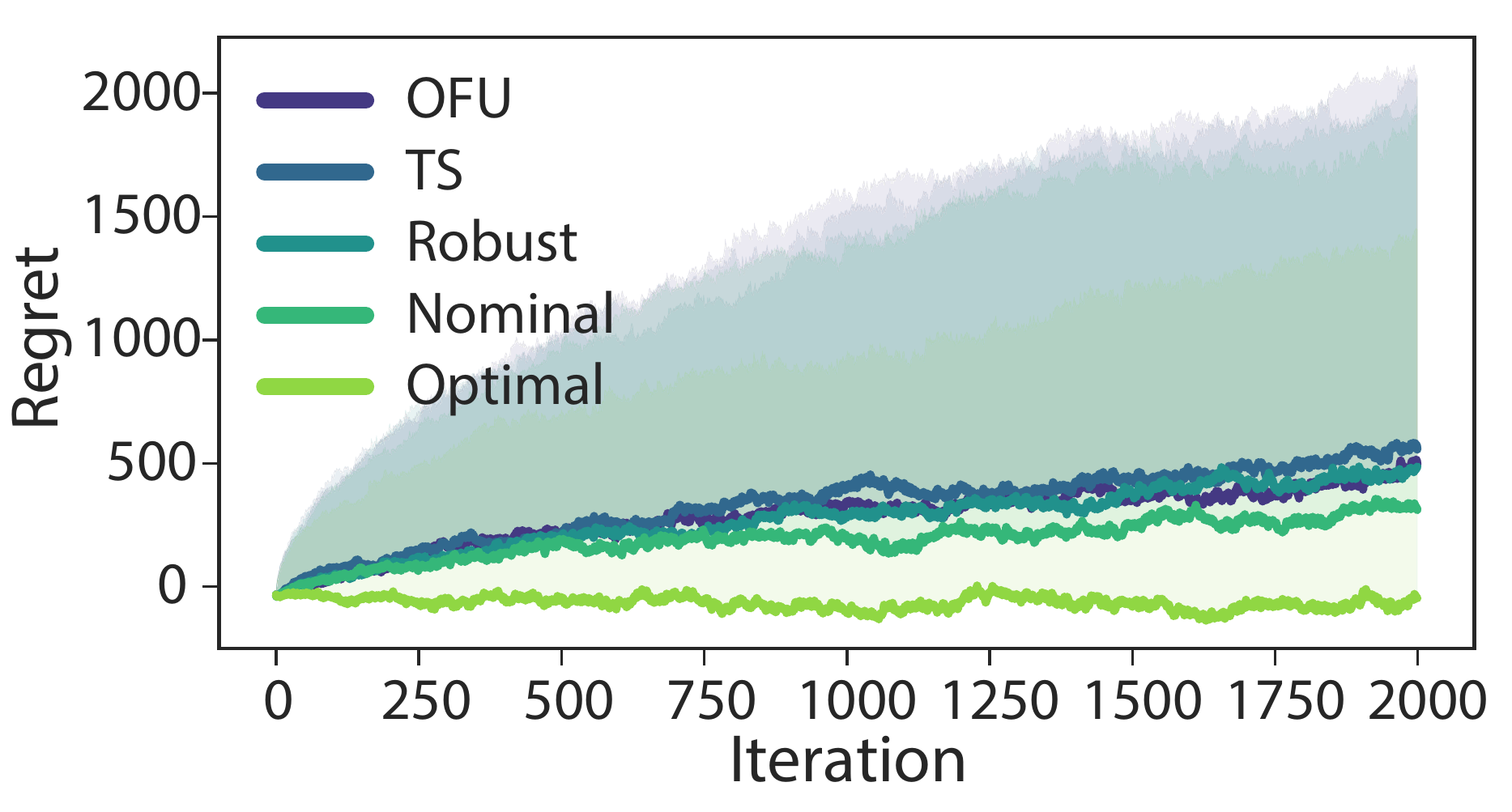}
\caption{A comparison of different adaptive methods on 500 experiments.  The median and 90th percentiles are shown.}
\label{fig:regrets}
\end{figure}
\end{example}


\section{Model-free Methods for LQR}
\label{sec:model_free}
We compare and contrast model-free methods from RL to
the model-based methods of the previous section.
We note that, while there is not a widely accepted technical definition
of a model-free method, we informally a method to be model-free
if it does not (as an intermediate step) estimate
the transition dynamics $(A, B)$. Instead, the model-free methods
we encounter either learn an alternative representation (such as
a value function), or directly search for the optimal controller.

\subsubsection{Overview of Model-free methods}

\paragraph{Approximate Policy Iteration}
\label{sec:model_free:lspi}

We start our discussion with a review of an important concept
from RL known as the \emph{state-action value function} or \emph{$Q$-function.}
Given a policy $\pi$, the (relative) state-action value function\footnote{We change the quadratic state matrix in the cost to $S$ here as not to confuse it with the $Q$ function.} is defined as
\begin{multline}
	Q^\pi(x, u) :=\\
	\lim_{T \to \infty} \E\left[ \sum_{t=0}^{T} (x_t^\T S x_t + u_t^\T R u_t - \lambda_\pi) \:\bigg|\: x_0 = x, u_0 = u \right] \:, \label{eq:relative_Q_function}
\end{multline}
where $\lambda_\pi$ is the infinite-horizon average cost of the policy $\pi$, achieved via $u_t = \pi(x_t)$.
The Bellman equation associated to \eqref{eq:relative_Q_function} is:
\begin{align}
	\lambda_\pi + Q^\pi(x, u) = c(x, u) + \E_{x' \sim p(\cdot|x, u)}[Q^\pi(x', \pi(x'))] \:. \label{eq:bellman_Q}
\end{align}
%
%
\emph{Policy iteration} (PI) is a classic algorithm from RL that works as follows.
Fix a starting policy $\pi_0$. Then, iteratively repeat: (1) Compute $Q_t$ as the state-value function for policy $\pi_t$,  (2) Update $\pi_{t+1}$ as $\pi_{t+1}(x) = \arg\min_{u} Q_t(x, u)$.
Step (1) of policy iteration generally requires knowledge of the transition dynamics.
Of course, $Q^\pi$ for a given policy $\pi$ can be estimated from data. One such method,
motivated by the Bellman equation \eqref{eq:bellman_Q}, is
based on temporal differences.
Choose a finite set of basis functions $\{ \phi_i \}_{i=1}^{K}$,
and suppose that $Q^\pi$ is well approximated in the span of $\phi_i$'s,
i.e. $Q^\pi \approx \sum_{i=1}^{K} w_i \phi_i$ for weights $w \in \R^K$.
Then given $T$ state transition tuples $\{ (x_i, u_i, x_{i+1}) \}_{i=1}^{T}$,
we can estimate the best fit weights $\wh$ as:
\begin{align}
	\wh = \left(  \sum_{t=1}^{T} \phi_t (\phi_t - \psi_{t+1})^\top \right)^{\dag} \sum_{t=1}^{T} \phi_t (c_t - \lambdah_\pi) \:.
\end{align}
Here, $\phi_t = \phi(x_t, u_t)$, where $\phi(x, u) \in \R^K$ is a column vector of each basis function
$\phi_i$ evaluated at the pair $(x, u)$. 
Furthermore, $\psi_t = \phi(x_t, \pi(x_t))$, $c_t$ is the instantaneous cost observed
for the $t$-th state transition, and $\lambdah_\pi$ is an estimate of $\lambda_\pi$, the infinite-horizon
average cost under the policy $\pi$.
This estimator $\wh$ is known as the \emph{least-squares temporal difference} estimator for $Q$-functions
(LSTD-Q) \cite{lagoudakis03}.
The LSTD-Q estimator is an \emph{off-policy} estimator, meaning that
the input $u_t$ applied to generate the data does \emph{not} have to follow the policy $\pi$.
This is one of the advantages of working with $Q$-functions (as opposed to value functions).
We now can state our first model-free algorithm, the least-squares policy iteration (LSPI)
algorithm of \cite{lagoudakis03}: (1) Estimate $\Qh_t \approx Q^{\pi_t}$ from data via LSTD-Q, (2) Update $\pi_{t+1}$ as $\pi_{t+1}(x) = \arg\min_{u} \Qh_t(x, u)$.

\paragraph{Derivative-free Search Methods}
\label{sec:model_free:random_search}

We now turn to methods which directly search for the optimal 
policy. Suppose we have a policy class $\Pi = \{\pi_\theta : \theta \in \Theta \}$
with $\Theta \subseteq \R^p$.
We are interested in solving the following optimization problem:
$\min_{\theta \in \Theta} J(\theta)$, where $J(\theta) = J(\pi_\theta)$ is the 
infinite-horizon average cost performance
of the policy $\pi_\theta$. In many problem formulations, the function
$\theta \mapsto J(\theta)$ is differentiable on $\Theta$. Therefore, in principle
one could run a local search method such as gradient descent:
$\theta_{t+1} = \theta_t - \eta \nabla_\theta J(\theta_t)$.
However, typically computing $\nabla_\theta J(\theta)$ requires knowledge of the dynamics.
To get around this, in many practical situations one can efficiently generate
(unbiased estimates) of $J(\theta)$ for any $\theta$ via rollouts or simulation.
Therefore, we can rely on a rich history of zero-th order optimization.
Here, we will focus on two popular zero-th order methods in RL.

\paragraph{Policy Gradients (REINFORCE)}

The policy gradient method was popularized by \cite{williams1992simple}, and
forms the foundation of many popular algorithms such as TRPO~\cite{schulman2015trust} and PPO~\cite{schulman17ppo}. The idea is to perturb the action sequence and use these perturbations to
estimate gradient information.
Let $J_\eta(\theta) = \lim_{T\to\infty} \E[\frac{1}{T}\sum_{t=1}^{T} c_t]$ where
$u_t = \pi_\theta(x_t) + \eta_t$ with $\eta_t \iid \calN(0, \sigma^2 I)$.
Now under sufficient regularity conditions that allow us to switch
the order of differentiation with the limit and the integral,
\begin{align*}
	\nabla_\theta J_\eta(\theta) &= \lim_{T\to\infty} \int_{\tau_{1:T}} \frac{1}{T} c(\tau_{1:T}) \nabla_\theta p(\tau_{1:T}) \: d\tau_{1:T} \\
	&=  \lim_{T\to\infty} \int_{\tau_{1:T}} \frac{1}{T} c(\tau_{1:T}) \nabla_\theta \log{p(\tau_{1:T})} p(\tau_{1:T}) \: d\tau_{1:T} \:.
\end{align*}
We now observe that because the transition dynamics are not a function of $\theta$:
$\nabla_\theta \log{p(\tau_{1:T})} = \sum_{t=1}^{T} \nabla_\theta \log{p(u_t|x_t)}$.
%
Furthermore, because $p(u_t|x_t) \stackrel{d}{=} \calN(\pi_\theta(x_t), \sigma^2 I)$,
we have that $\nabla_\theta \log{p(u_t|x_t)} = \frac{1}{\sigma^2} (D\pi_\theta(x_t))^\T \eta_t$,
where $D\pi_\theta(x_t) \in \R^{d \times p}$ is the Jacobian of the map $\theta \mapsto \pi_\theta(x_t)$.
Therefore:
\begin{align*}
	&\nabla_\theta J_\eta(\theta) = \lim_{T \to \infty} \E\left[ \frac{1}{T} \sum_{t=1}^{T} \frac{c(\tau_{1:T})}{\sigma^2} (D \pi_\theta(x_t))^\T \eta_t \right] \\
	&= \lim_{T \to \infty} \E\left[ \frac{1}{T} \sum_{t=1}^{T} \frac{c(\tau_{t:T})}{\sigma^2} (D \pi_\theta(x_t))^\T \eta_t \right] \\
	&= \lim_{T \to \infty} \E\left[ \frac{1}{T} \sum_{t=1}^{T} \frac{c(\tau_{t:T}) - b(\tau_{1:t-1}, x_t)}{\sigma^2} (D \pi_\theta(x_t))^\T \eta_t \right] \:,
\end{align*}
where the second and third equality both follow from iterating expectations and using 
the fact that $\E[\eta_t | x_1, \eta_1, ..., x_t] = 0$.
Here, $b(\tau_{1:t-1}, x_t)$ is a \emph{baseline} function which is chosen to reduce variance.
Therefore, we can form a gradient estimate of $J_\eta(\theta)$ by
choosing a large $T$, rolling out $T$ steps with $u_t = \pi_\theta(x_t) + \eta_t$,
and then forming $\widehat{g} = \frac{1}{T} \sum_{t=1}^{T} \frac{c(\tau_{t:T}) -  b(\tau_{1:t-1}, x_t)}{\sigma^2} (D \pi_\theta(x_t))^\T \eta_t$.

\paragraph{Random Search}

We now consider random finite differences.
There are many different variants of finite differences; we present
one of the simpler methods. The idea here is to perturb
the parameter space. Like for policy gradients, we consider
a surrogate function $J_\xi(\theta)$ defined as
$J_\xi(\theta) = \E_{\xi}[ J(\theta + \sigma \xi) ]$, where $\xi \sim \calN(0, I)$. 
It is a standard fact that the gradient $\nabla_\theta J_\xi(\theta)$ is
given by:
\begin{align*}
	\nabla_\theta J_\xi(\theta) = \E\left[ \frac{J(\theta + \sigma\xi) - J(\theta - \sigma\xi)}{2\sigma} \xi \right] \:.
\end{align*}
Hence we can construct a stochastic gradient by first choosing a large $T$, then sampling a random
perturbation $\xi$, and finally rolling out
a trajectory with $\pi_{\theta + \sigma\xi}$ and another with $\pi_{\theta - \sigma\xi}$,
and using the estimate $\widehat{g} = \frac{\frac{1}{T}\sum_{t=1}^{T} c_t - \frac{1}{T}\sum_{t=1}^{T} c'_t}{2\sigma}$.

\subsubsection{Experimental Evaluation}

We compare the previously described model-free methods to the model-based nominal control
(Section~\ref{sec:nominal}) as a baseline.
We consider the following LQR problem:
\begin{align*}
	A = \begin{bmatrix} 0.95 & 0.01 & 0 \\
	0.01 & 0.95 & 0.01 \\
	0 & 0.01 & 0.95
	\end{bmatrix}, \,
	B = \begin{bmatrix} 
	1 & 0.1 \\
	0 & 0.1 \\
	0 & 0.1 
	\end{bmatrix}, \,
	S = I_3,\, R = I_2.
\end{align*}
We choose an LQR problem where the $A$ matrix is stable, since
the model-free methods we consider require an initial
stabilizing controller; using a stable $A$ allows us to start
at $K_0 = 0_{2 \times 3}$.
We fix the process noise $\sigma_w = 1$.
As before, the model-based method learns $(A, B)$ using
least-squares,
exciting the system with Gaussian noise of variance $\sigma_u = 1$.

\begin{figure}[h!]
\centering
\includegraphics[width=0.5\columnwidth]{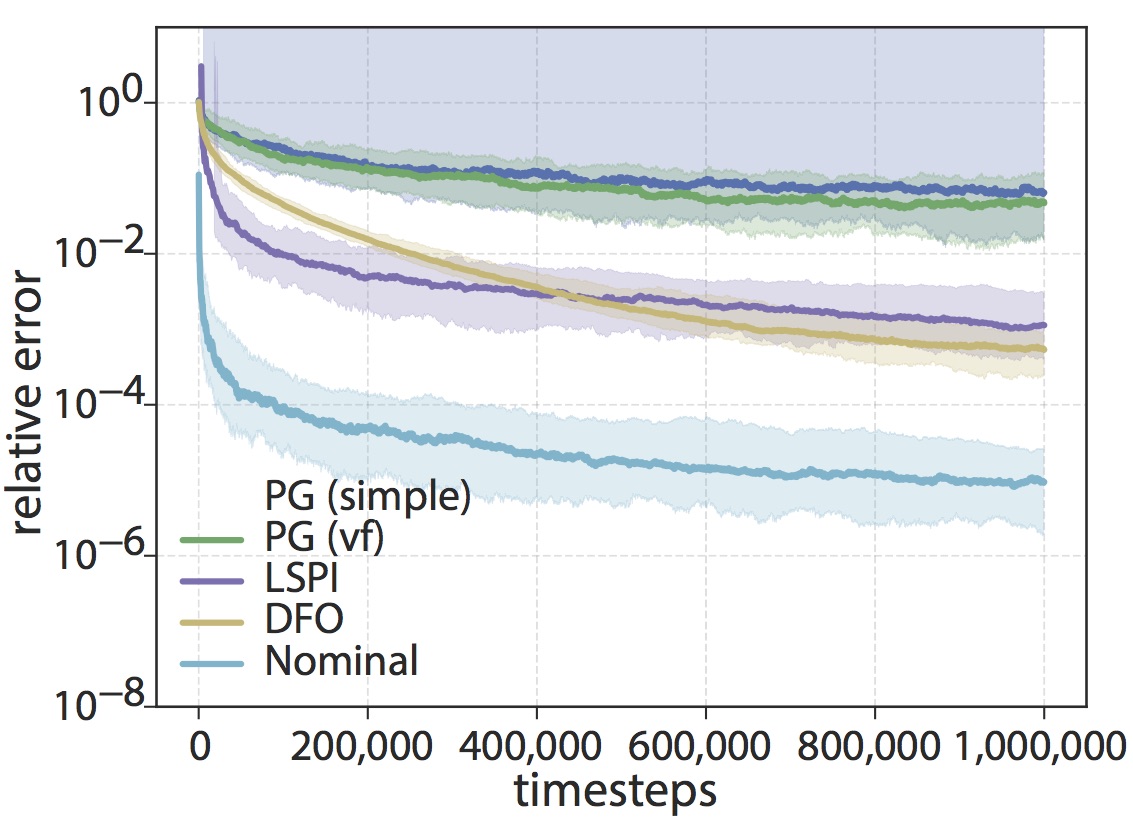}
\caption[Comparison of various model-free methods to nominal control.]{\small The performance of various model-free methods compared with
the nominal (Section~\ref{sec:nominal}) controller.
The shaded regions represent the lower 10th and upper 90th percentile over $100$ trials, and
the solid line represents the median performance.
Here, PG (simple) is policy gradients with the simple baseline, 
PG (vf) is policy gradients with the value function baseline,
PI is least-squares policy iteration,
and DFO is derivative-free optimization.
}
\label{fig:model_free:eval}
\end{figure}

For policy gradients and derivative-free optimization, we use the 
projected stochastic gradient descent (SGD) method with a constant step size $\mu$ as
the optimization procedure. The details of how we tuned each procedure
are described in~\cite{tu19thesis}.
For policy gradients we considered two different baselines. The 
\emph{simple baseline} is the one that uses the empirical average cost
$\frac{1}{T}\sum_{t=1}^{T} c_t$ of the previous iteration as the baseline.
The \emph{value function} baseline uses $b(x_t) = x_t^\top V x_t$ where
$V = \dlyap(A+BK, S+K^\top R K)$. Computing this $V$ requires
knowledge of the model $(A, B)$: a practical implementation would need to 
also estimate $V$, further degrading performance.

Figure~\ref{fig:model_free:eval} shows that the nominal model-based controller
substantially outperforms the model-free methods. We also see that the use of a baseline
reduces the variance for policy gradients, the necessity of which is
well understood in practice.

\subsubsection{Separation Results}

In the last section we saw that the model-based nominal method substantially
outperformed the model-free methods we considered in our experiment.
Can we make this rigorous? Here, we outline some results towards this,
based on \cite{tu2018gap}.
For these results, we consider a finite length $T$ horizon LQR problem with 
no input penalty:
$\min_{u_t} \E[\sum_{t=1}^{T} \norm{x_t}^2]$.
We also only consider dynamics which have the special property that
$\mathrm{range}(A) \subseteq \mathrm{range}(B)$ and $B$ has full column rank. These assumptions imply that the optimal solution 
is simply to cancel the state: $u_t = - B^{\dag} A x_t$. It means that
the optimal solution on a finite horizon is time-invariant,
which is typically not the case.

We first consider the risk of the model-based method.
The following theorem characterizes the performance of the
nominal method as the number of rollouts $N$ tends to infinity.
It is an asymptotic version of Theorem \ref{thm:ce_bound}
\begin{theorem}[Theorem 2.4, \cite{tu2018gap}]
\label{thm:model_free:CE}
The risk of the model-based nominal controller $\Kh$
on the $T$ length LQR problem described above satisfies:
\begin{align*}
	\lim_{N \to \infty} N \cdot \E[  J(\Kh) - J_\star ] = O(n_x n_u) + o_T(1) \:.
\end{align*}
\end{theorem}
Theorem~\ref{thm:model_free:CE} states that
$\E[ J(\Kh) - J_\star ] \approx O(nd/N)$ for large $N$.
Next, we study the behavior of the policy gradient algorithm.
We consider both a simple baseline $b(x_t) = \norm{x_t}^2$
and the value function baseline.
\begin{theorem}[Theorem 2.5, \cite{tu2018gap}]
\label{thm:model_free:PG}
The risk of the model-free policy gradient controller $\Kh$
on the $T$ length LQR problem described for the \emph{simple baseline}
satisfies:
\begin{align*}
	\liminf_{N \to \infty} N \cdot \E[  J(\Kh) - J_\star ] = \Omega(T^2 \cdot n_u n_x^3 )\:,
\end{align*}
and for the \emph{value function} baseline satisfies:
\begin{align*}
	\liminf_{N \to \infty} N \cdot \E[  J(\Kh) - J_\star ] = \Omega(T \cdot n_u n_x^2) \:.
\end{align*}
\end{theorem}
Above, the $O(\cdot)$ and $\Omega(\cdot)$ hides the dependence on
various properties of $(A, B)$ such as $\rho(A), \norm{B}_F, \sigma_{\min}(B)$ as well as
the noise variances $\sigma_w, \sigma_\eta$.
Comparing with Theorem~\ref{thm:model_free:CE},
Theorem~\ref{thm:model_free:PG} shows that policy gradient has
worse sample complexity than the model-based nominal controller:
$\E[ J(\Kh) - J_\star ] \approx \Omega(T^2 n_u n_x^3/N)$ for the simple baseline
and $\E[ J(\Kh) - J_\star ] \approx \Omega(T n_u n_x^2/N)$ for the value function baseline.
This result shows that, while more sophisticated baselines help, policy gradient
still suffers from worse sample complexity than the nominal method by 
factors of horizon length $T$ and state dimension $n$.

Finally, we turn to an information-theoretic lower bound which states that the model-based
nominal controller is optimal on this family of LQR instances.
\begin{theorem}[Theorem 2.6, \cite{tu2018gap}]
Consider the family of dynamics $\mathscr{G}(\rho, n_u) := \{ (\rho UU^\top, \rho U) \::\: U\in \R^{n_x\times n_u} \:, \:\: U^\top U = I \}$. Any algorithm $\mathcal{A}$ which plays
feedbacks of the form $u_t = K_i x_t + \eta_t$ with $\norm{K_i} \leq 1$ and $\eta_t \sim \calN(0, \sigma_u^2 I)$
incurs risk:
\begin{align*}
	\inf_{\mathcal{A}} \sup_{\rho \in (0, 1/4), (A,B) \in \mathscr{G}(d,\rho)} \E[J(\mathcal{A}) - J_\star] \gtrsim \frac{n_u(n_x-n_u)}{N} \:,
\end{align*}
if $n_u \leq n_x/2$ and $n_u(n_x-n_u)$ is greater than an absolute constant.
\end{theorem}
This result shows that the $O(n_xn_u/N)$ risk incurred by the nominal method is nearly optimal
up to constant factors over any algorithm which plays inputs of the form
$u_t = K x_t + \eta_t$, which includes both the nominal method, the policy gradient methods,
and the policy iteration method we described eariler.

\section{Conclusions}
\label{sec:conclusions}

This tutorial paper and our companion paper \cite{sysID} presented a broad overview of recent progress towards the finite-time analysis for reinforcement learning and self-tuning control methods.  We have attempted to provide a summary of representative results in this space that establish connections between the self-tuning control literature and methods recently proposed in reinforcement learning. The former are typically model-based, and are well-studied from a theoretical perspective, although more effort is still needed to better understand their finite-time behavior. The latter mostly adopt a model-free approach; they have had spectacular successes over the last few years, but lack strong theoretical guarantees, although researchers have been trying to validate design choices and algorithms a posteriori.  Empirically, there is rich evidence that learning algorithms exploiting prior knowledge about the system (such as a model parameterization) are more sample-efficient, but are also more sensitive to biases introduced by modeling errors.  Further, as we showed in Section \ref{sec:model_free}, there exists scenarios where there is a quantitative and provable gap between model-based and model-free methods.  A broader assessment of the advantages and drawbacks of modeling choices are however difficult to assess theoretically.  This survey summarized recent progress towards this goal within the limited scope of tabular MDPs and linear optimal control, but  we critically need to develop more broadly applicable tools towards the tighter analysis of learning algorithms. 

\bibliographystyle{IEEEtran}
{\bibliography{../learning-abridged.bib,../refsAP.bib}}

\end{document}